\documentclass[reqno,11pt]{amsart} 
\usepackage{amsmath}
\usepackage{amssymb}
\usepackage{amsthm}
\newcommand\NoBlackBoxes{\global\overfullrule0pt}
\NoBlackBoxes
\theoremstyle{plain} 

\def\4{\kern1pt}

\def\6{\vphantom0}

\def\8{\kern-10pt}
\def\7#1{_{(#1)}}

\makeatletter
\let\serieslogo@\relax
\let\@setcopyright\relax

\def\speciallabelmark#1{\def\@currentlabel{#1}}
\makeatother

\begin{document}

\def\ffrac#1#2{\raise.5pt\hbox{\small$\4\displaystyle\frac{\,#1\,}{\,#2\,}\4$}}
\def\ovln#1{\,{\overline{\!#1}}}
\def\ve{\varepsilon}
\def\kar{\beta_r}

\title{BERRY-ESSEEN BOUNDS IN THE ENTROPIC \\ 
CENTRAL LIMIT THEOREM
}

\author{S. G. Bobkov$^{1,4}$}
\thanks{1) School of Mathematics, University of Minnesota, USA;
Email: bobkov@math.umn.edu}
\address
{Sergey G. Bobkov \newline
School of Mathematics, University of Minnesota  \newline 
127 Vincent Hall, 206 Church St. S.E., Minneapolis, MN 55455 USA
\smallskip}
\email {bobkov@math.umn.edu} 

\author{G. P. Chistyakov$^{2,4}$}
\thanks{2) Faculty of Mathematics, University of Bielefeld, Germany;
Email: chistyak@math.uni-bielefeld.de}
\address
{Gennadiy P. Chistyakov\newline
Fakult\"at f\"ur Mathematik, Universit\"at Bielefeld\newline
Postfach 100131, 33501 Bielefeld, Germany}
\email {chistyak@math.uni-bielefeld.de}

\author{F. G\"otze$^{3,4}$}
\thanks{3) Faculty of Mathematics, University of Bielefeld, Germany;
Email: goetze@math.uni-bielefeld.de}
\thanks{4) Research partially supported by 
NSF grant DMS-1106530 and 
SFB 701}
\address
{Friedrich G\"otze\newline
Fakult\"at f\"ur Mathematik, Universit\"at Bielefeld\newline
Postfach 100131, 33501 Bielefeld, Germany}
\email {goetze@mathematik.uni-bielefeld.de}



\begin{abstract}
Berry-Esseen-type bounds for total variation and relative entropy
distances to the normal law are established  for the sums 
of non-i.i.d. random variables.
\end{abstract}

\maketitle
\markboth{S. G. Bobkov, G. P. Chistyakov and F. G\"otze}{Entropic Bounds}




\def\theequation{\thesection.\arabic{equation}}
\def\E{{\bf E}}
\def\R{{\bf R}}
\def\C{{\bf C}}
\def\P{{\bf P}}
\def\H{{\rm H}}
\def\Im{{\rm Im}}
\def\Tr{{\rm Tr}}

\def\k{{\kappa}}
\def\M{{\cal M}}
\def\Var{{\rm Var}}
\def\Ent{{\rm Ent}}
\def\O{{\rm Osc}_\mu}

\def\ep{\varepsilon}
\def\phi{\varphi}
\def\F{{\cal F}}
\def\L{{\cal L}}

\def\be{\begin{equation}}
\def\en{\end{equation}}
\def\bee{\begin{eqnarray*}}
\def\ene{\end{eqnarray*}}


\section{{\bf Introduction}}
\setcounter{equation}{0}

\vskip2mm
Let $X_1,\dots,X_n$ be independent (not necessarily identically distributed) 
random variables with mean $\E X_k = 0$ and finite variances 
$\sigma_k^2 = \E X_k^2$ $(\sigma_k > 0$). Put
$B_n = \sum_{k=1}^n \sigma_k^2$. 
Under additional moment assumptions, the normalized sum
$$
S_n = \frac{X_1 + \dots + X_n}{\sqrt{B_n}}
$$
has aproximately a standard normal distribution in a weak sense. 
Moreover, the closeness of the distribution function 
$F_n(x) = \P\{S_n \leq x\sqrt{B_n}\}$ to the standard 
normal distribution function
$$
\Phi(x) = \frac{1}{\sqrt{2\pi}} \int_{-\infty}^x e^{-y^2/2}\,dy
$$
has been studied intensively in terms of the so-called Lyapunov ratios
$$
L_s = \frac{\sum_{k=1}^n \E\, |X_k|^s}{B_n^{s/2}}.
$$
In particular, if all $X_k$ have finite third absolute moments, the 
classical Berry-Esseen theorem says that
\be
\sup_x |F_n(x) - \Phi(x)| \leq C L_3,
\en
where $C$ is an absolute constant (cf. e.g. [E], [F], [Pe]).

One of the most remarkable features of (1.1) is that the number
of summands does not {\it explicitly} appear in it, while
in the i.i.d. case, that is, when $X_k$ have equal distributions,
$L_3$ is of order $\frac{1}{\sqrt{n}}$, which is best possible
for the Kolmogorov distance under the 3-rd moment condition.

In this paper we shall prove bounds for stronger distances between 
$F_n$ and $\Phi$, such as total variation $\|F_n - \Phi\|_{{\rm TV}}$ 
and relative entropy $D(F_n||\Phi)$. However, these distances are
clearly useless for example when all summands have discrete 
distributions. Therefore, some further assumptions are needed.

When estimating the error of normal approximation by means of these
distances, it seems natural to require that every $X_k$ has an absolutely 
continuous distribution. Even with this assumption we cannot exclude
the case that our distances of $S_n$ to the normal law may be growing 
when the distributions of $X_k$ get near to discrete distributions.
Thus we shall assume that the densities of $X_k$ are bounded on 
a reasonably large part of the real line. This can be quaranteed
quite naturally, for instance, by using the entropy functional, 
defined for a random variable $X$ with density $p(x)$ by
$$
h(X) = -\int_{-\infty}^{+\infty} p(x) \log p(x)\,dx.
$$
Once $X$ has a finite second moment, the entropy is well-defined
as a Lebesgue integral, although the value $h(X) = -\infty$ is possible.
Introduce a related functional
$$
D(X) =
h(Z) - h(X) = \int_{-\infty}^{+\infty} p(x) 
\log \frac{p(x)}{\varphi_{a,\sigma}(x)}\,dx,
$$
where $Z$ is a normal random variable with density $\varphi_{a,\sigma}$ 
having the same mean $a$ and variance $\sigma^2$ as $X$. 
Note that this functional is affine invariant, that is,
$D(c_0 + c_1 X) = D(X)$, for all $c_0 \in \R$, $c_1 \neq 0$,
and in this sense it does not depend neither on the mean or the variance 
of $X$.

The quantity $D(X)$, denoted also $D(F_X||F_Z)$, where $F_X$ and $F_Z$ 
are the corresponding distributions of $X$ and $Z$, is known as the
"entropic distance to normality or Gaussianity". It may be characterized
as the shortest Kullback-Leibler distance from $F_X$
to the class of all normal laws on the real line.
In general, $0 \leq D(X) \leq +\infty$, and the equality $D(X) = 0$ is 
possible, when $X$ is normal, only. Moreover, by Pinsker's inequality,
the entropic distance dominates the total variation in the sense that
$$
D(X) \geq \frac{1}{2}\, \|F_X - F_Z\|_{{\rm TV}}^2.
$$

Thus, the size of $D(X)$ provides a strong distance of $F_X$ 
to normality, while finiteness of $D(X)$ guarantees that $F_X$ 
is separated from the class of discrete probability distributions.
Using $D$ for both purposes, one may obtain refinements of
Berry-Esseen's inequality (1.1) in terms of the total variation and
the entropic distances to normality for the distributions $F_n$.

\vskip5mm
{\bf Theorem 1.1.} {\it Let $D$ be a non-negative real number.
Assume that $X_k$ have finite third absolute moments, and
$D(X_k) \leq D$ $(1 \leq k \leq n)$. Then
\be
\|F_n - \Phi\|_{{\rm TV}} \leq C L_3,
\en
where the constant $C$ depends on $D$, only.
}

\vskip5mm
In particular, if all $X_k$ are equidistributed with $\E X_1^2 = 1$, 
we get
\be
\|F_n - \Phi\|_{{\rm TV}}\, \leq\, \frac{C}{\sqrt{n}}\, \E\, |X_1|^3
\en
with a constant $C$ depending on $D(X_1)$, only.
Although (1.2)-(1.3) seem to be new, related estimates in the 
i.i.d.-case were studied by many authors. For example, in the early 
1960's Mamatov and Sirazhdinov [M-S] found an exact asymptotic 
$\|F_n - \Phi\|_{{\rm TV}} = \frac{c}{\sqrt{n}} + o(\frac{1}{\sqrt{n}})$,
where the constant $c$ is proportional to $|\E X_1^3|$, and
which holds under the assumption that the distribution of $X_1$ has
a non-trivial absolutely continuous component (cf. also [Pr], [Se]).

Now, let us turn to the entropic distance to normality.

\vskip5mm
{\bf Theorem 1.2.} {\it Assume that $X_k$ have finite fourth absolute 
moments, and that $D(X_k) \leq D$ $(1 \leq k \leq n)$. Then
\be
D(S_n) \leq C L_4,
\en
where $C$ depends on $D$, only.
}

\vskip5mm
In (1.2) and (1.4) one may take $C = e^{c(D+1)}$, where
$c$ is an absolute constant. Moreover, $C$ can be chosen to be 
just a numerical constant, provided that $D$ 
is not too large, namely, if $D \leq c_0 \log \frac{1}{L_3}$ and 
$D \leq c_0 \log \frac{1}{L_4}$, respectively (with $c_0>0$ absolute).

These Berry-Esseen-type estimates are consistent
in view of the Pinsker-type inequality. In some sense, one may consider 
(1.4) as a stronger assertion than (1.2), which is indeed the case, when 
$L_4$ is of order $L_3^2$. (In general $L_3^2 \leq L_4$.)

In the i.i.d. case as in (1.3), the inequality (1.4) becomes
$$
D(S_n) \, \leq\, \frac{C}{n}\, \E X_1^4,
$$
where $C$ depends on $D(X_1)$ only. Thus, we obtain an error bound 
of order $O(1/n)$ under the 4th moment assumption. Note that 
the property $D(S_n) \rightarrow 0$ always holds under the second 
moment assumption (with finite entropy of $X_1$). This is the statement 
of the entropic central limit theorem, which is due to Barron [B]. Here, 
the convergence may have an arbitrarily slow rate. Nevertheless, 
the expected typical rate $D(S_n) = O(\frac{1}{n})$ was known 
to hold in some cases, for example, when $X_1$ has a distribution 
satisfying an integro-differential inequality of Poincar\'e-type. 
These results are due to Artstein, Ball, Barthe and Naor [A-B-B-N], 
and Barron and Johnson [B-J]; cf. also [J]. Recently, an exact
asymptotic for $D(S_n)$ has been studied in [B-C-G1]. If the entropy 
and the 4th moment of $X_1$ are finite, it was shown that
$$
D(S_n) = \frac{c}{n} + o\bigg(\frac{1}{n\log n}\bigg), \qquad 
c = \frac{1}{12}\,\big(\E X_1^3)^2.
$$
Moreover, with finite 3rd absolute moment
(and infinite 4th moment) such a relation may not hold, and it may 
happen that $D(S_n) \geq n^{-(1/2 + \ep)}$ for all $n$ large enough 
with a given prescribed $\ep>0$. This holds, for example, when
$X_1$ has density
$$
p(x) =  \int_0^{+\infty} \frac{1}{\sigma \sqrt{2\pi}}\, 
e^{-x^2/2\sigma^2}\,dP(\sigma),
$$
where $P$ is a probability measure on $(\frac{1}{e},+\infty)$ with density
$\frac{dP(\sigma)}{d\sigma} = (\sigma \log \sigma)^{-4}$
for $\sigma \geq e$ and with an arbitrary extension to the interval 
$\frac{1}{e} < \sigma < e$ satisfying
$\int_{1/e}^{+\infty} \sigma^2\,dP(\sigma) = 1$.

\vskip2mm
Therefore, in the general non-i.i.d.-case, the Lyapunov coefficient 
$L_3$ cannot be taken as an appropriate quantity for bounding the
error in Theorem 1.2, and $L_4$ seems more relevant. 
This is also suggested by the result of [A-B-B-N] for the weighted sums
$$
S_n = a_1 X_1 + \dots + a_n X_n \qquad (a_1^2 + \dots + a_n^2 = 1)
$$
of i.i.d. random variables $X_k$, such that 
$\E X_1 = 0$ and $\E X_1^2 = 1$. Namely, it is proved there that
\be
D(S_n)\, \leq\, \frac{L(a)}{c/2 + (1 - c/2)L(a)}\, D(X_1),
\en
where $L(a) = a_1^4 + \dots + a_n^4$ and $c \geq 0$ is an optimal 
constant in the Poincar\'e-type inequality 
$c\,\Var(u(X_1)) \leq \E\, [u'(X_1)]^2$.
But for the sequence $a_k X_k$ and $s=4$, the corresponding 
Lyapunov coefficient is exactly $L_4 = L(a)\, \E X_1^4$. 
Therefore, when $c = c(X_1)$ is positive, (1.5) yields the estimate
$$
D(S_n)\, \leq\, \frac{2D(X_1)}{c\,\E X_1^4}\,L_4,
$$
which is of a similar nature as (1.4).

Another interesting feature of (1.4) is that it may be connected
with transportation cost inequalities for the distributions $F_n$
of $S_n$ in terms of the quadratic Wasserstein distance $W_2$. For random
variables $X$ and $Z$ with finite second moments and distributions
$F_X$ and $F_Z$, this distance is defined by
$$
W_2^2(F_X,F_Z) = \inf_\pi \int_{-\infty}^{+\infty}\int_{-\infty}^{+\infty}
|x-y|^2\,d\pi(x,y),
$$
where the infimum is taken over all probability measures $\pi$ 
on the plane $\R^2$ with marginals $F_X$ and $F_Z$. The value 
$W_2^2(F_X,F_Z)$ is interpreted as the minimal expenses needed to
transport $F_Z$ to $F_X$, provided that it costs $|x-y|^2$ to move 
any "particle" $x$ to any "particle" $y$. 

The metric $W_2$ is of weak type in the sense that it can be used 
to metrize the weak convergence of probability distributions ([V]). 
Moreover, if $Z$ is standard normal and if $X$ has density, $W_2(F_X,F_Z)$ 
may be bounded in terms of the relative entropy by virtue of Talagrand's 
transportation inequality
\be
W_2^2(F_X,F_Z) \leq 2 D(F_X||F_Z)
\en
(cf. [T], or [B-G] for a different approach). If additionally $X$ 
has mean zero and unit variance, $D(F_X||F_Z) = D(X)$. Hence, applying
(1.6) with $X = S_n$, we get, by Theorem 1.2, 
\be
W_2(F_n,\Phi) \leq C \sqrt{L_4},
\en
where $C$ depends on $D$. In fact, this inequality
holds true with an absolute constant.
This result is due to Rio [Ri], who also studied more general
Wasserstein distances $W_r$, by relating them to Zolotarev's "ideal"
metrics. It has also been noticed in [Ri] that the 4-th moment condition
is essential, so the Laypunov's ratio $L_4$ in (1.7) cannot be replaced 
with $L_3$ including the i.i.d.-case (like in Theorem 1.2).

The paper is organized according to the following plan.

\vskip7mm
{\it Contents}

\vskip2mm
1. Introduction.

2. General bounds on total variation and entropic distance.

3. Entropic distance and Edgeworth-type approximation.

4. Quantile density decomposition.

5. Properties of the quantile decomposition.

6. Entropic bounds for Cramer constants of characteristic functions.

7. Repacking of summands.

8. Decomposition of convolutions.

9. Entropic approximation of $p_n$ by $\widetilde p_n$.

10. Integrability of characteristic functions $\widetilde f_n$
and their derivatives.

11. Proof of Theorem 1.1 and its refinement.

12. Proof of Theorem 1.2 and its refinement.

13. The case of bounded densities.

\vskip5mm


\vskip5mm
\section{{\bf General Bounds on Total Variation and Entropic Distance}}
\setcounter{equation}{0}

\vskip2mm
Let a random variable $X$ have an absolutely continuous
distribution $F$ with density $p(x)$ and finite first absolute moment.
We do not require that it has mean zero and/or unit variance.

First, we recall an elementary bound for the total variation distance 
$\|F - \Phi\|_{{\rm TV}}$ in terms of the characteristic function
$$
f(t) = \E\, e^{itX} = \int_{-\infty}^{+\infty} e^{itX}\,p(x)\,dx
\qquad (t \in \R).
$$
Introduce the characteristic function $g(t) = e^{-t^2/2}$ of the
standard normal law.

In the sequel, we use the notation
$$
\|u\|_2 = \bigg(\int_{-\infty}^{+\infty} |u(t)|^2\,dt\bigg)^{1/2}
$$
to denote the $L^2$-norm of a measurable complex-valued function $u$ 
on the real line (with respect to Lebesgue measure).

\vskip5mm
{\bf Proposition 2.1.} {\it We have
\be
\|F - \Phi\|_{{\rm TV}}^2 \, \leq \, \frac{1}{2}\,
\|f - g\|_2^2 + \frac{1}{2}\,\|f' - g'\|_2^2.
\en
}

This bound is standard (cf. e.g. [I-L], Lemma 1.3.1). 
In fact, the inequality (2.1) remains to hold for an arbitrary probability
distribution (in place of $\Phi$) with finite first absolute moment and characteristic function $g$. However, the general case
won't be needed in the sequel. 

Note that the assumption $\E\,|X| < +\infty$ guarantees that $f$ is 
continuously differentiable, so that the last integral in (2.1) 
makes sense.

Let $Z$ be a standard normal random variable, with density
$\varphi(x) = \frac{1}{\sqrt{2\pi}}\,e^{-x^2/2}$.
Consider the relative entropy
\be
D(X||Z) = D(F||\Phi) = \int_{-\infty}^{+\infty} p(x)\,\log
\frac{p(x)}{\varphi(x)}\,dx.
\en
As a preliminary bound, we first derive:

\vskip5mm
{\bf Lemma 2.2.} {\it For all $T \geq 0$,
\begin{eqnarray}
D(X||Z) 
 & \leq &
e^{-T^2/2} +
\sqrt{2\pi} \int_{-T}^T (p(x) - \varphi(x))^2\, e^{x^2/2}\ dx \nonumber \\
 & &
+ \ \frac{1}{2}\, \int_{|x| \geq T} x^2\,p(x)\, dx + 
\int_{|x| \geq T} p(x)\log p(x)\,dx. 
\end{eqnarray}
}

\vskip2mm
{\bf Proof.}
We split the integral in (2.2) into the two regions. 
For the interval $|x| \leq T$, using the elementary inequality
$t \log t \leq (t-1) + (t-1)^2$, $t \geq 0$, we have
\bee
\int_{-T}^T \frac{p}{\varphi} \log \frac{p}{\varphi}\ \varphi\,dx
 & \leq &
\int_{-T}^T \bigg(\frac{p}{\varphi} - 1\bigg)\, \varphi\,dx +
\int_{-T}^T \bigg(\frac{p}{\varphi} - 1\bigg)^2\, \varphi\,dx \\
 & = &
\int_{|x| \geq T} (\varphi - p)\, dx +
\int_{-T}^T \frac{(p - \varphi)^2}{\varphi}\ dx \\
 & = &
2\,(1 - \Phi(T)) - \int_{|x| \geq T} p(x)\, dx +
\sqrt{2\pi}
\int_{-T}^T (p(x) - \varphi(x))^2\, e^{x^2/2}\, dx.
\ene
For the second region, just write
\bee
\int_{|x| \geq T} p(x) \log \frac{p(x)}{\varphi(x)}\, dx
 & = &
\int_{|x| \geq T} p(x) \log p(x)\,dx \\
 & & + \ 
\log \sqrt{2\pi} \int_{|x| \geq T} p(x)\,dx + \frac{1}{2}
\int_{|x| \geq T} x^2\,p(x)\, dx.
\ene
It remains to collect these relations and use $\log \sqrt{2\pi} < 1$
together with a well-known elementary inequality
$1 - \Phi(T) \leq \frac{1}{2}\,e^{-T^2/2}$.
Thus, Lemma 2.2 is proved.

\vskip5mm
{\bf Remark.} If $p$ is bounded by a constant $M$, the estimate (2.3)
yields
\bee
D(X||Z) 
 & \leq &
e^{-T^2/2} +
\sqrt{2\pi} \int_{-T}^T (p(x) - \varphi(x))^2\, e^{x^2/2}\ dx \\
 & &
+ \ \frac{1}{2}\, \int_{|x| \geq T} x^2\,p(x)\, dx + 
\log M \int_{|x| \geq T} p(x)\,dx. 
\ene
This bound might be of interest in other applications,
although it involves the maximum of the density. 
For our purposes, the important integral in (2.3),
$
\int_{|x| \geq T} p(x)\log p(x)\,dx,
$
will be bounded in a different way and in terms of 
the characteristic functions, without involving the parameter $M$.


\vskip5mm
\section{{\bf Entropic Distance and Edgeworth-type Approximation}}
\setcounter{equation}{0}

\vskip2mm
To estimate the integrals in (2.3) in terms of the characteristic
functions like in Proposition 2.1, define
$$
\varphi_\alpha(x) = \varphi(x)\bigg(1 + \alpha\,\frac{x^3 - 3x}{3!}\bigg),
$$
where $\alpha$ is a parameter. These functions appear with $\alpha$
proportional to $n^{-1/2}$ in the Edgeworth-type expansions up 
to order 3 for densities of the normalized sums 
$S_n = \frac{X_1 + \dots + X_n}{\sqrt{B_n}}$ of i.i.d. summands. 
In the non-i.i.d. case such expansions hold as well
with 
$$
\alpha = \frac{1}{B_n^{3/2}}\, \sum_{k=1}^n \E X_k^3.
$$

Note that every $\varphi_\alpha$ has the Fourier transform
$$
g_\alpha(t) = \int_{-\infty}^{+\infty} e^{itx} \varphi_\alpha(x)\,dx =
g(t)\bigg(1 + \alpha\,\frac{(it)^3}{3!}\bigg),
$$
where $g(t) = e^{-t^2/2}$.

\vskip5mm
{\bf Proposition 3.1.} {\it Let $X$ be a random variable with
$\E\,|X|^3 < +\infty$. For all $\alpha \in \R$,
\be
D(X||Z) \ \leq \ \alpha^2 +
4\,\big(\|f - g_\alpha\|_2 + \|f''' - g_\alpha'''\|_2\big),
\en
where $Z$ is a standard normal random variable and
$f$ is the characteristic function of~$X$.
}

\vskip5mm
The assumption on the 3rd absolute moment is needed to insure
that $f$ has first three continuous derivatives.

As a particular case, the inequality (3.1) is valid for $\alpha = 0$,
as well. Then it becomes
$$
D(X||Z) \, \leq \, 4\,\big(\|f - g\|_2 + \|f''' - g'''\|_2\big),
$$
which may be viewed as a full analog of Proposition 2.1. However, with 
properly chosen values of $\alpha$, (3.1) may provide a much better 
asymptotic approximation (especially when applying it to the sums 
of independent random variables).

\vskip5mm
{\bf Proof.}
We may assume that the characteristic function $f$ and its first three 
derivatives are square integrable, so that the right-hand side
of (3.1) is finite. Note that in this case, $X$ has an absolutely
continuous distribution with some density $p$.

We apply Lemma 2.2. Given $T \geq 0$ to be specified later on, let us start 
with the estimation of the last integral in (2.3). Define the even function 
$\widetilde p(x) = p(x) + p(-x)$, so that 
$p \log p \leq p \log^+ \widetilde p$ (where we use the notation 
$a^+ = \max\{a,0\}$).
Subtracting $\varphi_\alpha(x)$ from $p(x)$ and then adding, one can write
\bee
\int_{|x| \geq T} p(x)\log p(x)\,dx
 & \leq &
\int_{|x| \geq T} p(x)\log^+ \widetilde p(x)\,dx \\
 & \leq &
\int_{-\infty}^{+\infty} |p(x) - \varphi_\alpha(x)|\, 
\log^+ \widetilde p(x)\,dx
 + \int_{|x| \geq T} \varphi_\alpha(x)\, \log^+ \widetilde p(x)\,dx.
\ene
But the function $\varphi_\alpha - \varphi$ is odd, so the last 
integral does not depend on $\alpha$ and is equal to
\be
\int_{|x| \geq T} \varphi(x)\, \log^+ \widetilde p(x)\,dx.
\en
To estimate it from above, one may use Cauchy's inequality together
with the elementary bound $(\log^+ t)^2 \leq Ct$, where
the optimal constant $C$ is equal to $4e^{-2}$. Since
$\int_{-\infty}^{+\infty} \widetilde p(x)\,dx = 2$,
(3.2) does not exceed
$$
\bigg(\int_{|x| \geq T} \varphi(x)^2\,dx\bigg)^{1/2}
\bigg(\int_{|x| \geq T} \big(\log^+ \widetilde p(x)\big)^2\,dx\bigg)^{1/2}
\leq
\bigg(\int_{|x| \geq T} \varphi(x)^2\,dx\bigg)^{1/2}\,
\frac{2\sqrt{2}}{e}.
$$
On the other hand,
$$
\bigg(\int_{|x| \geq T} \varphi(x)^2\,dx\bigg)^{1/2} =
\bigg(\frac{1}{\sqrt{\pi}}\,\left(1 - \Phi(T\sqrt{2})\right)\!\!\bigg)^{1/2}
\leq \, \frac{1}{\pi^{1/4}\sqrt{2}}\,e^{-T^2/2},
$$
where we applied the inequality
$1 - \Phi(x) \leq \frac{1}{2}\,e^{-x^2/2}$ ($x \geq 0$). Thus, 
using $\frac{2\sqrt{2}}{e} \cdot \frac{1}{\pi^{1/4}\sqrt{2}} < 1$
to simplify the constant, we get
$$
\int_{|x| \geq T} p(x)\log p(x)\,dx \leq
\int_{-\infty}^{+\infty} |p(x) - \varphi_\alpha(x)| \,
\log^+ \widetilde p(x)\,dx + e^{-T^2/2}.
$$

Here, again by the Cauchy inequality, the last integral does not exceed
$$
\frac{2\sqrt{2}}{e}
\bigg(\int_{-\infty}^{+\infty} (p(x) - \varphi_\alpha(x))^2\,dx\bigg)^{1/2}
=
\frac{2\sqrt{2}}{e} \cdot \frac{1}{\sqrt{2\pi}}
\bigg(\int_{-\infty}^{+\infty} |f(t) - g_\alpha(t)|^2\,dt\bigg)^{1/2},
$$
where we applied Plancherel's formula. The constant in front of the last
integral is smaller than $\frac{1}{2}$, so we arrive at the estimate
\be
\int_{|x| \geq T} p(x)\log p(x)\,dx \, \leq \,
\frac{1}{2}\, \|f - g_\alpha\|_2 + e^{-T^2/2}.
\en

Now, let us turn to the pre-last integral in (2.3). Once more,
subtracting $\varphi_\alpha(x)$ from $p(x)$ and then adding, one can write
$$
\int_{|x| \geq T} x^2 p(x)\,dx \leq \int_{-\infty}^{\infty} 
x^2\,|p(x) - \varphi_\alpha(x)|\,dx + 
\int_{|x| \geq T} x^2\,\varphi_\alpha(x)\,dx.
$$
Since the function $\varphi_\alpha - \varphi$ is odd, the last 
integral is equal to 
$$
\int_{|x| \geq T} x^2 \varphi(x)\,dx \, = \,
\frac{2}{\sqrt{2\pi}} \int_T^{+\infty} x^2\,e^{-x^2/2}\,dx \, = \,
2(1 - \Phi(T)) + \frac{2}{\sqrt{2\pi}}\ T e^{-T^2/2}
$$
(by direct integration by parts). Hence, using 
$2(1 - \Phi(T)) \leq e^{-T^2/2}$ once more, we get
\begin{eqnarray}
\frac{1}{2}\,\int_{|x| \geq T} x^2 p(x)\,dx
 & \leq &
\frac{1}{2}\,\int_{-\infty}^{+\infty} x^2\,|p(x) - \varphi_\alpha(x)|\,dx 
\nonumber \\
 & & 
+ \ \frac{1}{2}\, e^{-T^2/2} + \frac{1}{\sqrt{2\pi}}\ T e^{-T^2/2}.
\end{eqnarray}

In addition, by Cauchy's inequality, 
\bee
\bigg(\int_{-\infty}^{\infty} x^2\,|p(x) - \varphi_\alpha(x)|\,dx\bigg)^2
 & \leq &
\int_{-\infty}^{+\infty} \frac{dx}{1 + x^2}\,
\int_{-\infty}^{+\infty} (1+x^2)\,x^4\,(p(x) - \varphi_\alpha(x))^2\,dx \\
 & = &
\pi \int_{-\infty}^{+\infty} (x^4+x^6)\,(p(x) - \varphi_\alpha(x))^2\,dx \\
 & \leq &
\pi \int_{-\infty}^{+\infty} (1+2x^6)\,(p(x) - \varphi_\alpha(x))^2\,dx.
\ene
But, by Plancherel's formula,
\begin{eqnarray}
\int_{-\infty}^{+\infty} (p(x) - \varphi_\alpha(x))^2\,dx
 & = &
\frac{1}{2\pi}\,\|f - g_\alpha\|_2^2 \\
\int_{-\infty}^{+\infty} x^6\,(p(x) - \varphi_\alpha(x))^2\,dx
 & = &
\frac{1}{2\pi}\,\|f''' - g_\alpha'''\|_2^2.
\end{eqnarray}
Hence, 
\bee
\int_{-\infty}^{+\infty} x^2\,|p(x) - \varphi_\alpha(x)|\,dx 
 & \leq &
\bigg(\frac{1}{2}\,\|f - g_\alpha\|_2^2 + \|f''' - g_\alpha'''\|_2^2\bigg)^{1/2}
 \\
 & \leq &
\|f - g_\alpha\|_2 + \|f''' - g_\alpha'''\|_2,
\ene
and from (3.4),
\begin{eqnarray}
\frac{1}{2}\,\int_{|x| \geq T} x^2 p(x)\,dx
 & \leq &
\frac{1}{2}\, e^{-T^2/2} + \frac{1}{\sqrt{2\pi}}\ T e^{-T^2/2} \nonumber \\
 & & 
+ \ \frac{1}{2}\,\|f - g_\alpha\|_2 + \frac{1}{2}\,\|f''' - g_\alpha'''\|_2.
\end{eqnarray}
Using the bounds (3.3) and (3.7) in the inequality (2.3), we therefore 
obtain that
\begin{eqnarray}
D(X||Z) 
 & \leq &
\frac{5}{2}\, e^{-T^2/2} + \frac{1}{\sqrt{2\pi}}\ T e^{-T^2/2} \nonumber \\
 & & \hskip-15mm + \
\sqrt{2\pi} \int_{-T}^T (p(x) - \varphi(x))^2\, e^{x^2/2}\, dx +
\|f - g_\alpha\|_2 + \|f''' - g_\alpha'''\|_2.
\end{eqnarray}

Next, let us consider the integral in (3.8). First, writing
$$
p(x) - \varphi(x) = \big(p(x) - \varphi_\alpha(x)\big) + 
\alpha\,\frac{x^3 - 3x}{3!}\,\varphi(x)
$$
and applying an elementary inequality 
$(a+b)^2 \leq \frac{a^2}{1-t} + \frac{b^2}{t}$ ($a,b \in \R$, $0<t<1$)
with $t = 1/6$, we get
$$
(p(x) - \varphi(x))^2 \leq \frac{6}{5}\,
\big(p(x) - \varphi_\alpha(x)\big)^2 + \alpha^2\,
\frac{(x^3 - 3x)^2}{6}\,\varphi(x)^2,
$$
or equivalently,
$$
(p(x) - \varphi(x))^2 \, e^{x^2/2} \leq \frac{6}{5}\,
\big(p(x) - \varphi_\alpha(x)\big)^2 \, e^{x^2/2} +
\frac{1}{\sqrt{2\pi}}\ \alpha^2\,\frac{(x^3 - 3x)^2}{6}\,\varphi(x).
$$
Integrating this inequality over the interval $[-T,T]$ and
using $\E\, (Z^3 - 3Z)^2 = 6$, where $Z \sim N(0,1)$, we obtain
$$
\sqrt{2\pi} \int_{-T}^T (p(x) - \varphi(x))^2\, e^{x^2/2}\, dx 
\ \leq \ \frac{6}{5}\,\sqrt{2\pi} \int_{-T}^T 
(p(x) - \varphi_\alpha(x))^2\, e^{x^2/2}\, dx + \alpha^2.
$$

To estimate the last integral, first note that the function 
$t \rightarrow e^{t/2}/(2+t)$ is increasing for $t \geq 0$. Hence, 
for all $|x| \leq T$,
$$
e^{x^2/2} = \frac{e^{x^2/2}}{2 + x^2}\,\big(2 + x^2\big) 
\leq \frac{e^{T^2/2}}{2 + T^2}\,\big(3 + x^6\big),
$$
and thus, using (3.5)-(3.6),
\bee
\int_{-T}^T (p(x) - \varphi_\alpha(x))^2\, e^{x^2/2}\, dx
 & \leq & 
\frac{e^{T^2/2}}{2 + T^2} 
\int_{-T}^T (3 + x^6)\,(p(x) - \varphi_\alpha(x))^2\,dx \\
 & \leq & 
\frac{3}{2\pi}\ \frac{e^{T^2/2}}{2 + T^2}
\ \big(\|f - g_\alpha\|_2^2 + \|f''' - g_\alpha'''\|_2^2\big).
\ene
Putting $\ep = \|f - g_\alpha\|_2 + \|f''' - g_\alpha'''\|_2$, we get
$$
\sqrt{2\pi} \int_{-T}^T (p(x) - \varphi_\alpha(x))^2\, e^{x^2/2}\, dx \leq
\frac{18}{5\sqrt{2\pi}} \ \frac{e^{T^2/2}}{2 + T^2}\, \ep^2 + \alpha^2.
$$
Inserting this inequality in (3.8) leads to
\be
D(X||Z) \, \leq \,
\frac{5}{2}\, e^{-T^2/2} + \frac{1}{\sqrt{2\pi}}\ Te^{-T^2/2} +
\frac{18}{5\sqrt{2\pi}} \ \frac{e^{T^2/2}}{2 + T^2}\ \ep^2 + \ep + \alpha^2.
\en

It remains to optimize this bound over all $T \geq 0$. As before, consider
the function $\psi(t) = e^{t/2}/(2+t)$. It is increasing for $t \geq 0$
with $\psi(0) = \frac{1}{2}$. If $0 \leq \ep \leq 2$,
define $T = T_\ep$ to be the (unique) solution to the equation
$$
\psi(T^2) = \frac{1}{\ep}.
$$
In this case,
$$
T e^{-T^2/2} \cdot \frac{1}{\ep} = 
T e^{-T^2/2} \cdot \frac{e^{T^2/2}}{2 + T^2} \leq \frac{1}{2},
$$
so $T e^{-T^2/2} \leq \frac{\ep}{2}$. Furthemore, note that 
$$
e^{-T^2/2} \cdot \frac{1}{\ep} = e^{-T^2/2} \cdot
\frac{e^{T^2/2}}{2 + T^2}\leq \frac{1}{2},
$$
so $e^{-T^2/2} \leq \frac{\ep}{2}$. Applying these bounds in (3.9),
we arrive at
$$
D(X||Z) \, \leq \,
\frac{5\ep}{4} + \frac{1}{\sqrt{2\pi}}\, \frac{\ep}{2} +
\frac{18}{5\sqrt{2\pi}}\ \ep + \ep + \alpha^2 \, \leq \,
4\,\ep + \alpha^2,
$$
which is exactly the desired inequality (3.1).

In case $\ep \geq 2$, let us return to (3.8) and apply it with $T=0$.
This yields
$$
D(X||Z) \leq \frac{5}{2} + \ep < 4\,\ep, 
$$
which is even better than (3.1). Thus, Proposition 3.1 is proved.


\section{{\bf Quantile Density Decomposition}}
\setcounter{equation}{0}

\vskip2mm
In order to effectively apply Propositions 2.1 and 3.1, one has 
to solve two different tasks. The first one is to estimate 
integrals such as
$$
\int_{-T}^T |f(t) - g_\alpha(t)|^2\,dt, \quad
\int_{-T}^T |f'''(t) - g_\alpha'''(t)|^2\,dt
$$
over sufficiently large $t$-intervals with properly chosen values 
of the parameter $\alpha$. When the characteristic function $f$ has 
a multiplicative structure, i.e., corresponds to the sum of a large number 
of small independent summands, this task can be attacked by using 
classical Edgeworth-type expansions (for characteristic functions). 
Such expansions are well-known including the non-i.i.d. case, 
and we consider one of them in Section 12.

The second task concerns an estimation of integrals such as
$$
\int_{|x| \geq T} |f(t)|^2\,dt, \quad 
\int_{|x| \geq T} |f'''(t)|^2\,dt,
$$
which in general do not need to be small or even finite. The finiteness 
is quaranteed, for example, when $f$ is the Fourier transform of 
a bounded density $p$. For some purposes such as obtaining local 
limit theorems, it is therefore natural to restrict oneself to the 
case of bounded densities. For other purposes, such as an estimation 
of the total variation or relative entropy, the density $p$ may slightly 
be modified, so that the new density, say $\widetilde p$, will be bounded,
and at the same time will only slightly change the total variation
distance or relative entropy with respect to the standard normal law.

To this aim, we shall use the so-called quantile density 
decomposition, based on the following elementary observation.
(In fact, it is needed in case of bounded densities, as well.)

\vskip5mm
{\bf Proposition 4.1.} {\it Let $X$ be a random variable with density 
$p$. Given $0 < \kappa < 1$, the real line can be partitioned into 
two Borel sets $A_0, A_1$ such that $p(x) \leq p(y)$, for all 
$x \in A_0$, $y \in A_1$, and
$$
\int_{A_0} p(x)\,dx = \k, \qquad
\int_{A_1} p(x)\,dx = 1-\k.
$$
}

The argument is based on the continuity of the measure $p(x)\,dx$ and 
is omitted.

Clearly, for some real number $m_\k$ we get
$$
A_0 \subset \{x \in \R: p(x) \leq m_\k\}, \qquad
A_1 \subset \{x \in \R: p(x) \geq m_\k\}.
$$
Here, $m_\k$ represents a quantile (or one of the quantiles) for the 
function $p$ viewed as a random variable on the probability space 
$(\R,p(x)\,dx)$. In other words, $m_\k = m_\k(p(X))$ is a quantile 
of order $\k$ for the random variable $p(X)$. If $\k = \frac{1}{2}$, 
the index is usually omitted, and then $m = m(p(X))$ denotes 
a median of $p(X)$.

\vskip5mm
{\bf Definition 4.2.}
Define the densities $p_0$ and $p_1$ to be the normalized restrictions of 
$p$ to the sets $A_0$ and $A_1$, respectively. As a result, 
we have an equality
\be
p(x) = \k p_0(x) + (1-\k)\, p_1(x),
\en
which we call the quantile density decomposition for $p$
(respectively -- the median density decomposition, when $\k = \frac{1}{2}$).

\vskip5mm
Let us mention one obvious, but important property of the functionals 
$m_\k(p(X))$, assuming that $X$ has a finite second moment.

\vskip5mm
{\bf Proposition 4.3.} {\it The functionals 
$$
Q_\k(X) = m_\k(p(X)) \sqrt{\Var(X)}
$$ 
are affine invariant. That is, for all $a \in \R$ and $b \neq 0$,
$Q_\k(a + bX) = Q_\k(X)$.
}

\vskip5mm
More precisely, one should either assume in the latter equality
that the quantile $m_\k(p(X))$ is determined uniquely, or to use
specific quantiles satisfying the relation
$m_\k(p_{a,b}(a + bX)) = |b|^{-1}\, m_\k(p(X))$, where 
$p_{a,b}$ denotes the density of the random variable $a + bX$.


\section{{\bf Properties of the Quantile Decomposition}}
\setcounter{equation}{0}

\vskip2mm
In this section we establish basic properties of the quantile density 
decomposition. Although for purposes of Theorems 1.1-1.2
the median decomposition is sufficient, the general case is no more
difficult (but may be used to provide more freedom especially 
for improving $D$-dependent constants).

First, let us bound from above the quantiles $m_\k = m_\k(p(X))$ 
in terms of the entropic distance to normality. 

\vskip5mm
{\bf Proposition 5.1.} {\it Let $X$ be a random variable with finite 
variance $\sigma^2$ $(\sigma>0)$, having an absolutely continuous 
distribution, and let $0 < \kappa < 1$. Then
$$
m_\kappa \leq \frac{1}{\sigma\sqrt{2\pi}}\ e^{(D(X) + 1)/(1 - \kappa)}.
$$
In particular,
$$
m \leq \frac{1}{\sigma\sqrt{2\pi}}\,e^{2D(X) + 2}.
$$
}

\vskip2mm
{\bf Proof.} By Proposition 4.3, we may assume that $X$ has mean zero and 
variance one. Let $A = \{x \in \R: p(x) \geq m_\kappa\}$. 
By the definition of the quantiles,
$$
\int_A p(x)\,dx \geq 1-\kappa.
$$
Since $p(x) \geq m_\kappa$ on the set $A$, we have
\bee
\int_{-\infty}^{+\infty} p(x)\log\bigg(1 + \frac{p(x)}{\varphi(x)}\bigg)\,dx 
 & \geq &
\int_A p(x)\log\bigg(1 + \frac{m_\kappa}{\varphi(x)}\bigg)\,dx \\
 & \geq &
\int_A p(x)\log\frac{m_\kappa}{\varphi(x)}\,dx \\
 & = &
\log(m_\kappa\sqrt{2\pi}) \int_{A} p(x)\,dx + \frac{1}{2} \int_A x^2 p(x)\,dx \\
 & \geq &
(1-\kappa)\, \log(m_\kappa\sqrt{2\pi}).
\ene
On the other hand, using an elementary inequality 
$t \log(1+t) - t \log t \leq 1$ ($t \geq 0$), we get
\bee
\int_{-\infty}^{+\infty} p(x)\log\bigg(1 + \frac{p(x)}{\varphi(x)}\bigg)\,dx 
 & = &
\int_{-\infty}^{+\infty} \frac{p(x)}{\varphi(x)}\log\bigg(1 + \frac{p(x)}{\varphi(x)}\bigg)\,\varphi(x)\,dx \\
 & \leq &
\int_{-\infty}^{+\infty} \frac{p(x)}{\varphi(x)}\log \frac{p(x)}{\varphi(x)} \ 
\varphi(x)\,dx + 1
 \ = \
D(X) + 1. 
\ene 
Hence, $(1-\kappa) \log(m_\kappa\sqrt{2\pi}) \leq D(X) + 1$, 
and the proposition follows.

\vskip5mm
Now, let $V_0$ and $V_1$ be random variables with densities $p_0$ and 
$p_1$ from the quantile decomposition (4.1). They have means 
$a_j = \E\,V_j$ and variances $\sigma_j^2 = \Var(V_j)$, connected by
$$
\k a_0 + (1-\k)\, a_1 = \E X,
$$
and
\be
\left(\k a_0^2 + (1-\k)\,a_1^2\right) + 
\left(\k \sigma_0^2 + (1-\k)\,\sigma_1^2\right) = \E X^2,
\en
provided that $X$ has a finite second moment.

The next step is to prove upper bounds for the entropies of 
$V_0$ and $V_1$.

\vskip5mm
{\bf Proposition 5.2.} {\it If $X$ has mean zero and finite second moment, 
then
$$
\k D(V_0) + (1-\k)\,D(V_1) \leq D(X) -\k\log \k - (1-\k) \log(1-\k).
$$
In particular, in case of the median decomposition,
$$
D(V_0) + D(V_1) \leq 2D(X) + 2\log 2.
$$
}

\vskip2mm
{\bf Proof.} Let $\Var(X) = \sigma^2$ $(\sigma > 0$). We may assume that 
$D(X)$ is finite. By Definition 4.2,
\bee
-h(V_0)
 & = & 
\int_{-\infty}^{+\infty} p_0(x) \log p_0(x)\, dx \\
 & = &
\int_{A_0} (p(x)/\k) \log(p(x)/\k)\,dx \ = \
-\log \k + \frac{1}{\k} \int_{A_0} p(x) \log p(x)\,dx,\\
\ene
and similarly, 
$-h(V_1) = -\log(1-\k) + \frac{1}{1-\k} \int_{A_1} p(x) \log p(x)\,dx$. 
Adding the two equalities with weights, we get
\be
-\k h(V_0) - (1-\k)\,h(V_1) =  -\k\log \k - (1-\k) \log(1-\k) - h(X).
\en
Recall that
\bee
D(V_0) & = & h(Z_0) - h(V_0), \qquad {\rm where} \ \ Z_0 \sim N(a_0,\sigma_0^2), \\
D(V_1) & = & h(Z_1) - h(V_1), \qquad {\rm where} \ \ Z_1 \sim N(a_1,\sigma_1^2), \\
D(X) & = & h(Z) \, - \, h(X), \qquad \, {\rm where} \ \ Z \sim N(0,\sigma^2).
\ene
Hence, from (5.2),
\bee
\k D(V_0) + (1-\k)\,D(V_1) & = & 
\k h(Z_0) + (1-\k)\,h(Z_1) \\
 & & 
-\k\log \k - (1-\k) \log(1-\k) + (D(X) - h(Z)) \\
 & = & 
\k \log(\sigma_0 \sqrt{2\pi e}\,) + 
(1-\k)\,\log(\sigma_1 \sqrt{2\pi e}\,) \\
 & & 
-\k\log \k - (1-\k) \log(1-\k) + 
\big(D(X) - \log(\sigma\sqrt{2\pi e}\,)\big) \\
 & = & 
-\k\log \k - (1-\k) \log(1-\k) + D(X) + 
\log \frac{\sigma_0^\k \sigma_1^{1-\k}}{\sigma}.
\ene

Finally, by (5.1), and the arithmetic-geometric inequality,
$$
\sigma_0^{2\k} \sigma_1^{2(1-\k)} \leq \k \sigma_0^2 + (1-\k)\,\sigma_1^2 
\leq \sigma^2,
$$
so, $\frac{\sigma_0^\k \sigma_1^{1-\k}}{\sigma} \leq 1$.
Proposition 5.2 is proved.

\vskip5mm
Note that bounds on $D(X)$ provide a quantitative measure of
non-degeneracy of the distributions of $V_j$ via positivity of their
variances $\sigma_j^2$.

\vskip5mm
{\bf Proposition 5.3.} {\it Let $X$ be a random variable with mean zero 
and variance $\sigma^2$ $(\sigma>0)$, having finite entropy. Then
$$
\sigma_0 > \sigma\, e^{-(D(X) + 4)/\k}, \qquad
\sigma_1 > \sigma\, e^{-(D(X) + 4)/(1-\k)}.
$$
}

{\bf Proof.} By homogeneity with respect to $\sigma$, one may assume that 
$\sigma = 1$. 

We modify the argument from the proof of Proposition 5.1. First note that
\begin{eqnarray}
\log(\sigma_0 \sqrt{2\pi e}\,) & = &
D(V_0) - \int_{-\infty}^{+\infty} p_0(x) \log p_0(x)\,dx \nonumber \\
 & \geq &
- \int_{-\infty}^{+\infty} p_0(x) \log p_0(x)\,dx
 \ = \
- \int_{A_0} (p(x)/\k)\, \log(p(x)/\k)\,dx \nonumber \\
 & = &
\log \k - \frac{1}{\k} \int_{A_0} p(x)\, \log p(x)\ dx,
\end{eqnarray}
where $A_0$ is a set from Definition 4.2.

In order to estimate the last integral, put $r(x) = e^{-a^2 x^2/2}$ 
with parameter $a>0$. Using the property $r(x) \leq 1$ and once more 
the inequality $t\log(1+t) \leq t\log t + 1$ $(t \geq 0)$, we get
\bee
\int_{A_0} p(x)\, \log p(x)\ dx 
 & \leq &
\int_{-\infty}^{+\infty} p(x)\, \log\bigg(1 + \frac{p(x)}{r(x)}\bigg)\ dx \\
 & = &
\int_{-\infty}^{+\infty} \frac{p(x)}{r(x)}\, 
\log\bigg(1 + \frac{p(x)}{r(x)}\bigg)\ r(x)\,dx \\ 
 & \leq &
\int_{-\infty}^{+\infty} \bigg[\frac{p(x)}{r(x)}\, 
\log \frac{p(x)}{r(x)} + 1\bigg]\ r(x)\,dx \\ 
 & = &
\int_{-\infty}^{+\infty} p(x)\,\log p(x)\,dx + \frac{a^2}{2} 
\int_{-\infty}^{+\infty} p(x)\,x^2\, dx + \int_{-\infty}^{+\infty} r(x)\,dx \\
 & = &
D(X) - \log\big(\sqrt{2\pi e}\,\big) + \bigg(\frac{a^2}{2} + \frac{1}{a}\, \sqrt{2\pi}\bigg).
\ene
The right-hand side is minimized for $a = (2\pi)^{1/6}$ in which case
we obtain that
$$
\int_{A_0} p(x)\, \log p(x)\ dx  \leq 
D(X) - \log(\sqrt{2\pi e}\,) + \frac{3}{2}\,(2\pi)^{1/3} <
D(X) + 1.35.
$$
Together with (5.3), the above estimate yields
$$
\log(\sigma_0 \sqrt{2\pi e}\,) > \log \k - \frac{1}{\k}\, \big(D(X) + 1.35\big).
$$
But $\log(\sqrt{2\pi e}\,) \sim 1.42 < \frac{1.42}{\k}$, so
$\log \sigma_0 > \log \k - \frac{1}{\k}\, (D(X) + 2.77)$, or equivalently,
$$
\sigma_0 > \k\, e^{-(D(X) + 2.77)/\k}.
$$
Finally, using $\k > e^{-1/\k}$, the above estimate may be simplified to
$$
\sigma_0 > e^{-(D(X) + 3.77)/\k},
$$
which gives the first estimate on $\sigma_0$.
The second estimate for $\sigma_1$ is similar.

Thus, Proposition 5.3 is proved. Note that in case of the median 
decomposition, it becomes
$$
\sigma_0 > c \sigma\, e^{-2D(X)}, \qquad
\sigma_1 > c \sigma\, e^{-2D(X)},
$$
where $c$ is a positive absolute constant. One may take
$c = e^{-8}$, for example.


\vskip5mm
\section{{\bf Entropic Bounds for Cramer constants 
of Characteristic Functions}}
\setcounter{equation}{0}

\vskip2mm
If a random variable $X$ has an absolutely continuous distribution
with density, say $p$, then, by the Riemann-Lebesgue theorem,
its characteristic function
$$
f(t) = \E\, e^{itX} = \int_{-\infty}^{+\infty}e^{itx} p(x)\,dx \qquad
(t \in \R)
$$
satisfies
$f(t) \rightarrow 0$, as $t \rightarrow \infty$. Hence, for all $T>0$,
$$
\delta_X(T) = \sup_{|t| \geq T} |f(t)| < 1.
$$

An important problem is how to quantify this separation property 
(that is, separation from 1) by giving explicit upper bounds on 
the quantity $\delta_X(T)$, sometimes called Cramer constant.
(At least $\delta_X(T) < 1$ is refered to as Cramer's condition (C)).
This problem arises naturally in local limit theorems for densities of 
the sums of non-identically distributed independent summands.
Furthermore, it appears in the study of bounds and rates of 
convergence in the central limit theorem for strong metrics 
including the total variation and relative entropy. 
For our purposes, it is desirable to bound $\delta_X(T)$ explicitly 
in terms of the entropy of $X$ or, what is more relevant, in terms of 
the entropic distance to normality $D(X)$. Thus, this quantity 
controls separation of the distribution of $X$ from the class 
of discrete measures on the line.

A preliminary answer may be given in terms of the variance 
$\sigma^2 = \Var(X)$, when it is finite, and in cases
where the density $p$ is uniformly bounded.

\vskip5mm
{\bf Proposition 6.1.} {\it Assume $p(x) \leq M$ a.e. Then,
for all $t$ real,
\be
|f(t)|\,\leq\, 1 - c\,\frac{\min\{1,\sigma^2 t^2\}}{M^2 \sigma^2},
\en
where $c > 0$ is an absolute constant.
}

\vskip5mm
In a slightly different form, this bound was obtained in the mid 1960's by 
Statulevi\v{c}ius [St]. He also considered more complicated quantities
reflecting the behavior of the density $p$ on non-overlapping intervals 
of the real line.

The inequality (6.1) can be generalized by involving
non-bounded densities, but then $M$ should be replaced by other 
quantites such as quantiles $m_\kappa = m_\kappa(p(X))$
of the random variable $p(X)$. One can also remove any assumption 
on the moments of $X$ by replacing the standard deviation by the quantiles 
of the random variable $X-X'$, where $X'$ is an independent copy of $X$.
We refer to [B-C-G2] for details, where the following bound is derived.

\vskip5mm
{\bf Proposition 6.2.} {\it Let $X$ be a random variable with finite
variance $\sigma^2$ and finite entropy. Then, for all $t$ real,
\be
|f(t)|\,\leq\, 1 - c\,\min\{1,\sigma^2 t^2\}\, e^{-4D(X)},
\en
where $c > 0$ is an absolute constant.
}

\vskip5mm
At the expense of a worse constant in the exponent, this bound can be derived 
directly from (6.1) by combining it with Propositions 5.1 and 5.3.

Indeed, we may assume that $\E X = 0$.
Let $V_0$ and $V_1$ be random variables with densities 
$p_0$ and $p_1$ from the median decomposition (4.1), that is, for
$\kappa = \frac{1}{2}$, and
denote by $f_0$ and $f_1$ the corresponding characteristic functions, 
so that $f = \frac{1}{2}\, f_0 + \frac{1}{2}\,f_1$. Hence, for all $t$,
\be
|f(t)| \, \leq \, \frac{1}{2}\, |f_0(t)| + \frac{1}{2}.
\en
Since $p_0$ is bounded -- more precisely, $p_0(x) \leq m = m(p(X))$, 
one can apply Proposition 6.1 to the random variable $V_0$ with $M = m$.
Then (6.1) and (6.3) give
$$
|f(t)|\,\leq\, 
1 - c\,\frac{\min\{1,\sigma_0^2 t^2\}}{m^2 \sigma_0^2},
$$
where $\sigma_0^2 = \Var(V_0)$ and $c>0$ is an absolute constant.
Note that $\sigma_0^2 \leq 2\sigma^2$, according to (5.1).

Now, by Proposition 5.1,
$$
m^2 \sigma_0^2 \leq 2 m^2 \sigma^2 \leq \frac{1}{\pi}\, e^{4 D(X) + 4}.
$$
Hence,
$$
|f(t)|\,\leq\, 1 - c_1\,\min\{1,\sigma_0^2 t^2\}\, e^{-4D(X)}.
$$
Finally, by Propositions 5.3, $\sigma_0^2 > c_2 \sigma^2\, e^{-4D(X)}$, 
so
$$
|f(t)|\,\leq\, 
1 - c_3\,\min\{1,\sigma^2 t^2\}\, e^{-8\,D(X)}
$$
with some absolute constants $c_j>0$.


\vskip5mm
\section{{\bf Repacking of Summands}}
\setcounter{equation}{0}

\vskip2mm
We now consider a sequence of independent (not necessarily identically 
distributed) random variables $X_1,\dots,X_n$ and their sum
$S_n = X_1 + \dots + X_n$. Let $\E X_k = 0$, $\E X_k^2 = \sigma_k^2$ 
($\sigma_k > 0$). One may always assume without loss of 
generality that $\sigma_1^2 + \dots + \sigma_n^2 = 1$, 
so that $\Var(S_n) = 1$.

In addition, all $X_k$ are assumed to have absolutely continuous
distributions, having finite entropies in each place, where the
functional $D$ is used.

To study integrability properties of the characteristic function $f_n$ 
of $S_n$ (more precisely -- of its slightly modified variants 
$\widetilde f_n$), it will be more convenient to work with 
a different representation,
$$
S_n = V_1 + \dots + V_N,
$$
where the new independent summands represent appropriate partial sums
of the $X_l$ resulting in almost equal variances, such that at the
same time the number of blocks, $N$, is still reasonably large. 
Such a representation may be introduced just by taking
\be
V_k \ = \sum_{n_{k-1} < l \leq n_k} X_l,
\en
where $n_0 = 0$ and
$n_k = \max\{\,l \leq n:\, \sigma_1^2 + \dots + \sigma_l^2 \leq \frac{k}{N}\}$.

The number of new summands is restricted in terms of the parameter
$$
\sigma = \max_l \sigma_l
$$
which in general may be an arbitrary real number between 
$\frac{1}{\sqrt{n}}$ and 1.

\vskip5mm
{\bf Lemma 7.1.} {\it If $N \leq \frac{1}{2 \sigma^2}$, then for each 
$k = 1,\dots,N$,
\be
\frac{1}{2N} < \Var(V_k) < \frac{2}{N}.
\en
}

{\bf Proof.} If $n_1 = n$, then necessarily $N = 1$ and $V_1 = S_n$,
so (7.2) holds immediately. 

If $n_1 < n$, then, by the definition, 
$\Var(V_1) \leq \frac{1}{N}$ and $\Var(V_1 + X_{n_1 + 1}) > \frac{1}{N}$.
The latter implies 
$\Var(V_1) > \frac{1}{N} - \sigma^2 \geq \frac{1}{2N}$,
thus proving (7.2) for $k=1$.

Now, let $2 \leq k \leq N$. Again by the definition, 
$\Var(S_{n_k}) \leq \frac{k}{N}$ and
$\Var(S_{n_{k-1} + 1}) > \frac{k-1}{N}$. The latter implies
$\Var(S_{n_{k-1}}) > \frac{k-1}{N} - \sigma^2$.
Combining the two bounds, we get
$$
\Var(V_k) = \Var(S_{n_k}) - \Var(S_{n_{k-1}}) \leq 
\frac{k}{N} - \bigg(\frac{k-1}{N} - \sigma^2\bigg) = \frac{1}{N} + \sigma^2 
< \frac{2}{N}.
$$
On the other hand,
$$
\Var(V_k) > \bigg(\frac{k}{N} - \sigma^2\bigg) - \frac{k-1}{N} 
= \frac{1}{N} - \sigma^2 \geq \frac{1}{2N}.
$$
Lemma 7.1 is proved.

\vskip2mm
Thus, to obtain the property (7.2), it seems suggestive to take 
$N = [\frac{1}{2 \sigma^2}]$ (the integer part). However, this choice
is not used in the proof of Theorems 1.1-1.2, since we need
to express $N$ as a suitable function of Lyapunov's coefficients.

As another useful property of the representation (7.1), let us mention
the following.

\vskip5mm
{\bf Lemma 7.2.} {\it If $\max_{l \leq n} D(X_l) \leq D$, then
$\max_{k \leq N} D(V_k) \leq D$, as well.
}

\vskip5mm
This is due to the general bound $D(X+Y) \leq \max\{D(X),D(Y)\}$,
which holds for arbitrary independent random variables with finite 
second moments and absolutely continuous distributions. It can easily
be derived, for example, from the entropy power inequality
$$
e^{2 h(X + Y)} \geq e^{2 h(X)} + e^{2 h(Y)},
$$
cf. [C-D-T].

Now, let $\rho_k$ denote density of the random variable $V_k$.
For each $\rho_k$, one may consider a median density decomposition
\be
\rho_k(x) = \frac{1}{2}\, \rho_{k0}(x) + \frac{1}{2}\,\rho_{k1}(x)
\en
in accordance with Definition 4.2 for the parameter $\kappa = \frac{1}{2}$.

In particular, $\rho_{k0}(x) \leq m$, where $m = m(\rho_k(V_k))$
is a median of the random variable $\rho_k(V_k)$. Note that by 
Proposition 5.1 with $X = V_k$ and Lemmas 7.1-7.2, if 
$\max_{j \leq n} D(X_j) \leq D$, we immediately obtain that
\be
m(\rho_k(V_k)) \, \leq \, \frac{1}{v_k\sqrt{2\pi}} \, e^{2D + 2} \, \leq \,
\sqrt{N} \, e^{2D + 2},
\en
where $v_k = \sqrt{\Var(V_k)}$.

Let $V_{kj}$ be random variables with densities $\rho_{kj}$ and
characteristic functions
$$
\hat \rho_{kj}(t) = \E\, e^{it V_{kj}} = 
\int_{-\infty}^{+\infty} e^{itx}\,\rho_{kj}(x)\,dx, \quad j = 0,1.
$$
We collect their basic properties in the following lemma.

\vskip5mm
{\bf Lemma 7.3.} {\it Assume that $N \leq \frac{1}{2 \sigma^2}$\, and\,
$\max_{l \leq n} D(X_l) \leq D$.  For all $k \leq N$ and $j = 0,1$,

\vskip4mm
$a)$\ $D(V_{kj}) \leq 2D + 2$,

\vskip2mm
$b)$\ $\Var(V_{kj}) > \frac{1}{2N}\, e^{-4(D+4)}$,

\vskip2mm
$c)$\ $|\hat \rho_{kj}(t)| \leq 1 - c\,e^{-12\,D}$ for all 
$|t| \geq \sqrt{N}$ with an absolute constant $c>0$.
}

\vskip5mm
{\bf Proof.} The first assertion follows from Lemma 7.2 and
Proposition 5.2 applied with $X = V_k$. For the second one, combine
Proposition 5.3 with $X = V_k$ and Lemmas 7.1-7.2 to get
$$
v_{kj} > v_k\, e^{-2(D(V_k)+4)} \geq v_k\, e^{-2(D+4)}
\geq \frac{1}{\sqrt{2N}}\, e^{-2(D+4)},
$$
where $v_{kj}^2 = \Var(V_{kj})$ ($v_{kj}>0$).
For the assertion in $c)$, combine Proposition 6.2 for $X = V_{kj}$ and 
the previous steps, which give
\bee
|\hat \rho_{kj}(t)|
 & \leq & 
1 - c\,\min\big\{1,v_{kj}^2 t^2\big\}\, e^{-4D(V_{kj})} \\
 & \leq & 
1 - c\,\min\big\{1,t^2/(2N)\big\}\, e^{-4(D+4)} e^{-4 (2D+2)} \\
 & \leq & 
1 - c'\,\min\big\{1,t^2/N\big\}\, e^{-12\,D}
\ene
with some absolute constants $c,c' > 0$.


\vskip5mm
\section{{\bf Decomposition of Convolutions}}
\setcounter{equation}{0}

\vskip2mm
Starting from the representation $S_n = V_1 + \dots +V_N$
with the summands defined in (7.1), one can write the density of
$S_N$ as the convolution
$$
p_n = \rho_1 * \dots * \rho_N,
$$
where $\rho_k$ denotes the density of $V_k$. Moreover, a direct 
application of the median decomposition (7.3) leads to the representation
$$
p_n = 2^{-N}
\sum\,(\rho_{10}^{\delta_1} * \rho_{11}^{1-\delta_1}) * \dots * 
(\rho_{N0}^{\delta_N} * \rho_{N1}^{1-\delta_N}),
$$
where the summation is carried out over all $2^N$ sequences $\delta_k$ 
with values 0 and 1.

Let an integer number $m_0 \geq 0$ be given (For our purposes, one may take 
$m_0 = 3$). For $N \geq m_0 + 1$, we split the above sum into the two parts, 
so that
$$
p_n = q_{n0} + q_{n1},
$$ 
where
$$
q_{n0} \ = 2^{-N} \
\sum_{\delta_1 + \dots + \delta_N > m_0}
(\rho_{10}^{\delta_1} * \rho_{11}^{1-\delta_1}) * \dots * 
(\rho_{N0}^{\delta_N} * \rho_{N1}^{1-\delta_N}),
$$
$$
q_{n1} \ = 2^{-N} \
\sum_{\delta_1 + \dots + \delta_N \leq m_0} 
(\rho_{10}^{\delta_1} * \rho_{11}^{1-\delta_1}) * \dots * 
(\rho_{N0}^{\delta_N} * \rho_{N1}^{1-\delta_N}).
$$
Put
$$
\ep_n \, = \, \int_{-\infty}^{+\infty} q_{n1}(x)\,dx 
 \, = \, 2^{-N} \sum_{k=0}^{m_0}\, \frac{N!}{k!\, (N-k)!}.
$$
One can easily see that
\be
\ep_n \, \leq \, 2^{-(N-1)}\,N^{m_0}.
\en

\vskip2mm
{\bf Definition 8.1.} Put
\be
\widetilde p_n(x) = p_{n0}(x) = \frac{1}{1 - \ep_n}\, q_{n0}(x),
\en
and similarly $p_{n1}(x) = \frac{1}{\ep_n}\,q_{n1}(x)$. Thus, we get 
the decomposition
\be
p_n(x) = (1 - \ep_n) p_{n0}(x) + \ep_n p_{n1}(x).
\en
Accordingly, introduce the associated characteristic functions
$$
\widetilde f_n(t) = f_{n0}(t) = 
\int_{-\infty}^{+\infty} e^{itx}\widetilde p_n(x)\,dx,
\qquad f_{n1}(t) = \int_{-\infty}^{+\infty} e^{itx} p_{n0}(x)\,dx.
$$

\vskip2mm
The probability densities $\widetilde p_n(x) = p_{n0}(x)$ are bounded 
and provide a strong approximation for $p_n(x)$. 
Indeed, from (8.3) it follows that
\be
|\widetilde p_n(x) - p_n(x)| = \ep_n |p_{n0}(x) - p_{n1}(x)|
\en
which together with the bound (8.1) immediately implies:

\vskip5mm
{\bf Proposition 8.2.} {\it For all $n \geq N \geq m_0 + 1$,
$$
\int_{-\infty}^{+\infty} |\widetilde p_n(x) - p_n(x)|\,dx \, \leq \,
2^{-(N-2)}\,N^{m_0}.
$$
In particular, the corresponding characteristic functions satisfy,
for all $t \in \R$,
$$
|\widetilde f_n(t) - f_n(t)|\, \leq \, 2^{-(N-2)}\,N^{m_0}.
$$
}

We need a similar inequality for derivatives of characteristic 
functions. To this aim, we shall use absolute moments 
$\E\, |X_k|^s$ and the associated Lyapunov ratios
$$
L_s = \sum_{k=1}^n \E\, |X_k|^s \quad (s \geq 2).
$$ 

Let $V_{kj}$ ($1 \leq k \leq N$, $j=0,1$) be independent random 
variables with respective densities $\rho_{kj}$ from the median
decomposition (7.3) for the random variables $V_k$. For each sequence 
$\delta = (\delta_k)_{1 \leq k \leq N}$ 
with values 0 and 1, the convolution
$$
\rho^{(\delta)} =
(\rho_{10}^{\delta_1} * \rho_{11}^{1-\delta_1}) * \dots * 
(\rho_{N0}^{\delta_N} * \rho_{N1}^{1-\delta_N})
$$
represents the density of the sum
$$
S(\delta) = \sum_{k=1}^N\, \delta_k V_{k0} + (1-\delta_k) V_{k1}.
$$
If all moments $\E\, |X_k|^s$ are finite, (7.3) yields
\be
\E\, |V_k|^s = \frac{1}{2}\,\E\, |V_{k0}|^s + \frac{1}{2}\,\E\, |V_{k1}|^s.
\en
Hence, for the $L^s$-norm
$\|S(\delta)\|_s = \big(\E\, |S(\delta)|^s\big)^{1/s}$, using the
Minkowski inequality, we have
\bee
\|S(\delta)\|_s
 & \leq &
\sum_{k=1}^N \|\delta_k V_{k0} + (1-\delta_k) V_{k1}\|_s \\
 & \leq &
\sum_{k=1}^N \big(\delta_k\, \|V_{k0}\|_s + (1-\delta_k)\, \|V_{k1}\|_s\big)
 \ \leq \
2^{1/s} \sum_{k=1}^N \|V_k\|_s,
\ene
where (8.5) was used in the last step. But
$$
\frac{1}{N}\,\sum_{k=1}^N \|V_k\|_s =
\frac{1}{N}\,\sum_{k=1}^N \big(\E\, |V_k|^s\big)^{1/s} \leq 
\bigg(\frac{1}{N}\,\sum_{k=1}^N \E\, |V_k|^s\bigg)^{1/s},
$$
so
$$
\E\, |S(\delta)|^s \leq 2 N^{s-1} \sum_{k=1}^N \E\, |V_k|^s
\leq 2 N^s\, \E\, |S_n|^s,
$$
where we used $\E\, |V_k|^s \leq \E\, |S_n|^s$ 
(due to Jensen's inequality). 

\vskip2mm
Write
$\E\, |S(\delta)|^s = \int_{-\infty}^{+\infty} |x|^s\, 
\rho^{(\delta)}(x)\,dx$.
Recalling the definition of $q_{nj}$ and $\ep_n$, we get
\bee
\int_{-\infty}^{+\infty} |x|^s\, q_{n0}(x)\,dx
 & = &
2^{-N} \sum_{\delta_1 + \dots + \delta_N > m_0} \E\, |S(\delta)|^s
 \ \leq \ 
2\, \E\, |S_n|^s\, (1-\ep_n)\,N^s, \\
\int_{-\infty}^{+\infty} |x|^s\, q_{n1}(x)\,dx
 & = &
2^{-N} \sum_{\delta_1 + \dots + \delta_N \leq m_0} \E\, |S(\delta)|^s
 \ \leq \
2\, \E\, |S_n|^s\, \ep_n N^s.
\ene
Hence, by the definition of $p_{n0}$,
$$
\int_{-\infty}^{+\infty} |x|^s\, p_{n0}(x)\,dx \, \leq \, 
2\, \E\, |S_n|^s\, N^s,
$$
and similarly for $p_{n1}$. But, from (8.4),
$$
|x|^s\, |\widetilde p_n(x) - p_n(x)| \leq 
\ep_n |x|^s\, (p_{n0}(x) + p_{n1}(x)),
$$
so, applying (8.1),
$$
\int_{-\infty}^{+\infty} |x|^s\, |\widetilde p_n(x) - p_n(x)|\,dx \, \leq \, 
\E\, |S_n|^s\, 2^{-(N - 3)} N^{m_0 + s}.
$$

On the other hand, Rosenthal's inequality (cf. e.g. [Ro], [P-U]) gives
\be
\E\,|S_n|^s \leq C_s \bigg(1 + \sum_{j=1}^n \E\, |X_j|^s\bigg) =
C_s\,(1 + L_s), \qquad s \geq 2,
\en
with some constants $C_s$, depending on $s$, only (where the assumption 
$\E S_n^2 = 1$ is used). Note that in case $1 \leq s \leq 2$, 
there is also an obvious bound $\E\,|S_n|^s \leq 1$. 

One may summarize, using the constant $C_s$ in Rosenthal's inequality (8.6).

\vskip5mm
{\bf Proposition 8.3.} {\it Assume that $L_s$ is finite $(s \geq 2)$. 
For all $n \geq N \geq m_0 + 1$,
$$
\int_{-\infty}^{+\infty} |x|^s\, |\widetilde p_n(x) - p_n(x)|\,dx \, 
\leq \, C_s (1 + L_s)\, 2^{-(N - 3)}\, N^{m_0 + s}.
$$
In particular, if $s$ is an integer, the $s$-th derivative of the 
corresponding characteristic functions satisfies, for all $t$ real,
$$
|\widetilde f_n^{(s)}(t) - f_n^{(s)}(t)|\, \leq \,
C_s (1 + L_s)\,2^{-(N-3)}\,N^{m_0 + s}.
$$
}

\vskip2mm
For $s = 1$ and $s = 2$, it is better to use $\E\, |S_n| \leq 1$ 
and $\E S_n^2 = 1$ instead of (8.6). For $s=3$, Rosenthal's inequality 
can be shown to hold with constant $C_3 = 2$. Hence, we obtain:

\vskip5mm
{\bf Corollary 8.4.} {\it For all $n \geq N \geq m_0 + 1$ and $t \in \R$,
$$
|\widetilde f_n^{(s)}(t) - f_n^{(s)}(t)|\, \leq \, 
2^{-(N-3)}\, N^{m_0 + s} \qquad (s =1,2).
$$
Moreover, if $L_3$ is finite,
$$
|\widetilde f_n'''(t) - f_n'''(t)|\, \leq \, 
(1 + L_3)\,2^{-(N-4)}\, N^{m_0 + 3}.
$$
}


\vskip5mm
\section{{\bf Entropic Approximation of $p_n$ by $\widetilde p_n$}}
\setcounter{equation}{0}

\vskip2mm
As before, let $X_1,\dots,X_n$ be independent random variables 
with $\E X_k = 0$, $\E X_k^2 = \sigma_k^2$ ($\sigma_k>0$), such that 
$\sigma_1^2 + \dots + \sigma_n^2 = 1$. Moreover, let $X_k$ have
absolutely continuous distributions with finite entropies, and let
$p_n$ denote the density of the sum 
$$
S_n = X_1 + \dots + X_n.
$$ 
Put $\sigma^2 = \max_k \sigma_k^2$.

The next step is to extend the assertion of Propositions 8.2-8.3 
to relative entropies, with respect to the standard normal 
distribution on the real line with density 
$$
\varphi(x) = \frac{1}{\sqrt{2\pi}}\, e^{-x^2/2}.
$$
Thus put
$$
D_n = \int p_n(x) \log \frac{p_n(x)}{\varphi(x)}\ dx, \qquad
\widetilde D_n = \int \widetilde p_n(x) 
\log \frac{\widetilde p_n(x)}{\varphi(x)}\ dx.
$$

Recall that the modified densities $\widetilde p_n$ are constructed 
in Definition 8.1 with arbitrary integers $0 \leq m_0 < N \leq n$ 
on the basis of the representation (7.1), based on the independent 
random variables $V_k$ and the median decomposition (7.3) 
for the densities $\rho_k$ of $V_k$.

\vskip5mm
{\bf Proposition 9.1.} {\it Let $D = \max_{k} D(X_k)$. Given that
$m_0 + 1 \leq N \leq \frac{1}{2\sigma^2}$, we have
\be
|\widetilde D_n - D_n| < 2^{-(N-6)}\, N^{m_0+1} \,(D + 1).
\en
}

We shall use a few elementary properties of the convex function 
$L(u) = u \log u$ ($u \geq 0$).

\vskip5mm
{\bf Lemma 9.2.} {\it For all $u,v \geq 0$ and $0 \leq \ep \leq 1$,

\vskip4mm
$a)$ \ $L((1 - \ep)\,u + \ep v) \leq (1-\ep)\, L(u) + \ep L(v)$;

\vskip1mm
$b)$ \ $L((1 - \ep)\,u + \ep v) \geq (1-\ep)\, L(u) + \ep L(v) + 
u L(1-\ep) + v L(\ep)$.
}

\vskip5mm
{\bf Proof of Proposition 9.1.} Define
$$
D_{nj} = \int p_{nj}(x) \log \frac{p_{nj}(x)}{\varphi(x)}\ dx \qquad (j=0,1),
$$
so that $\widetilde D_n = D_{n0}$, where the densities $p_{nj}$ 
have been defined in (8.2)-(8.3).

By Lemma 9.2 $a)$,
$D_n \leq (1 - \ep_n)D_{n0} + \ep_n D_{n1}$. On the other hand, 
by Lemma 9.2 $b)$,
$$
D_n \geq \big((1 - \ep_n)D_{n0} + \ep_n D_{n1}\big) +
\ep_n \log \ep_n + (1-\ep_n) \log(1-\ep_n).
$$
The two estimates give
\be
|\widetilde D_n - D_n| \leq 
\ep_n (D_{n0} + D_{n1}) - \ep_n \log \ep_n - (1-\ep_n) \log(1-\ep_n).
\en
Hence, we need to give appropriate bounds on both $D_{n0}$ and $D_{n1}$.

To this aim, as before, let $V_{kj}$ ($1 \leq k \leq N$, $j=0,1$) be 
independent random variables with respective densities $\rho_{kj}$ 
from the median decomposition (7.3) for $V_k$.
By Definition 4.2, we have the identity (5.1), which for $V_k$
reads
$$
v_k^2 = \bigg(\frac{1}{2}\,a_{k0}^2 + \frac{1}{2}\,a_{k1}^2\bigg) + 
\bigg(\frac{1}{2}\,v_{k0}^2 + \frac{1}{2}\,v_{k1}^2\bigg),
$$
where $a_{kj} = \E V_{kj}$, $v_{kj}^2 = \Var(V_{kj})$ and
$v_k^2 = \Var(V_k)$. Using Lemma 7.1, this implies
\be
v_{k0}^2 \leq 2v_k^2 \leq \frac{4}{N}, \qquad  
v_{k1}^2 \leq 2v_k^2 \leq \frac{4}{N}.
\en

As in the previous section, for each sequence
$\delta = (\delta_k)_{1 \leq k \leq N}$ with values 0 and 1, consider the convolution
$$
\rho^{(\delta)} =
(\rho_{10}^{\delta_1} * \rho_{11}^{1-\delta_1}) * \dots * 
(\rho_{N0}^{\delta_N} * \rho_{N1}^{1-\delta_N}),
$$
i.e., the densities of the random variables
$$
S(\delta) = \sum_{k=1}^N\, \delta_k V_{k0} + (1-\delta_k) V_{k1}.
$$
By convexity of the function $u\log u$,
\begin{eqnarray}
\hskip-17mm D_{n1} \ \leq \
\frac{1}{\ep_n} \ 2^{-N} \sum_{\delta_1 + \dots + \delta_N \leq m_0}
\int_{-\infty}^{+\infty} \rho^{(\delta)}(x)\, \log \frac{\rho^{(\delta)}(x)}{\varphi(x)}\,dx, \\
\hskip-5mm D_{n0} \ \leq \
\frac{1}{1-\ep_n} \ 2^{-N} \sum_{\delta_1 + \dots + \delta_N > m_0}
\int_{-\infty}^{+\infty} \rho^{(\delta)}(x)\, 
\log \frac{\rho^{(\delta)}(x)}{\varphi(x)}\,dx.
\end{eqnarray}

Furthermore, if $S$ denotes a random variable with variance $v^2$ 
$(v>0)$ having density $\rho$, and if $Z$ is a standard normal random
variable, the relative entropy of $S$ with respect to $Z$ is connected 
with the entropic distance to normality $D(S)$ by the simple formula
\be
D(S||Z) = \int \rho(x) \log \frac{\rho(x)}{\varphi(x)}\ dx =
D(S) + \log\frac{1}{v} + \frac{\E S^2 - 1}{2}.
\en
In the case $S = S(\delta)$, applying Lemma 7.3 $b)$, we have 
$$
v^2 \, = \,
\sum_{k=1}^N \left[\,\delta_k v_{k0}^2 + (1-\delta_k)\, v_{k1}^2\right]
 \, \geq \,
\frac{1}{2}\, e^{-4(D + 4)},
$$
hence
\be
\log \frac{1}{v} \, \leq \, 2D + 9.
\en
In addition, arguing as in the proof of Proposition 8.2, specialized 
to the particular case $s=2$, and applying (9.3), we get
\bee
\|S(\delta)\|_2
 & \leq &
\sum_{k=1}^N \|\delta_k V_{k0} + (1-\delta_k) V_{k1}\|_2 \\
 & \leq &
\sum_{k=1}^N \big(\delta_k\, \|V_{k0}\|_2 + (1-\delta_k)\, \|V_{k1}\|_2\big)
 \ \leq \
\sqrt{2}\, \sum_{k=1}^N v_k \ \leq \ 2\sqrt{N}.
\ene
Hence, $\E S(\delta)^2 \leq 4N$. 
Combining this estimate with (9.7), we get that
$$
\log\frac{1}{v} + \frac{\E S(\delta)^2 - 1}{2} \, \leq \,
(2D + 9) + 2N.
$$

Consequently, if we apply this bound in (9.6) with $S = S(\delta)$, we obtain
\be
D(S(\delta)||Z) \leq D(S(\delta)) + (2D + 9) + 2N.
\en

The remaining term, $D(S(\delta))$, can be estimated 
by virtue of the same general inequality
$D(X+Y) \leq \max\{D(X),D(Y)\}$ mentioned before. This bound 
can be applied to all summands of $S(\delta)$, which together 
with Lemma 7.3 $a)$ gives
$$
D(S(\delta)) \leq \max_{1 \leq k \leq N} \max\{D(V_{k0}),D(V_{k1})\}
\leq 2D + 2.
$$
Applying this result in (9.8), we arrive at
$$
\int_{-\infty}^{+\infty} 
\rho^{(\delta)}(x)\, \log \frac{\rho^{(\delta)}(x)}{\varphi(x)}\,dx =
D(S(\delta)||Z) \leq 4D + 11 + 2N.
$$
Finally, by (9.4)-(9.5), we have similar bounds for $D_{n0}$ and $D_{n1}$,
namely,
$$
D_{n0} \leq 4D + 11 + 2N, \qquad
D_{n1} \leq 4D + 11 + 2N.
$$

Having obtained these estimates, we are prepared to return to (9.2), 
which thus gives
\be
|\widetilde D_n - D_n| \, \leq \, 2\ep_n\, (4D + 11 + 2N) + 
\ep_n \log\frac{1}{\ep_n} + (1-\ep_n) \log\frac{1}{1-\ep_n}.
\en

To simplify this bound, consider the function
$H(\ep) = \ep \log\frac{1}{\ep} + (1-\ep) \log\frac{1}{1-\ep}$, which
is defined for $0 \leq \ep \leq 1$, is concave and symmetric about the
point $\frac{1}{2}$, where it attains its maximum $H(\frac{1}{2}) = \log 2$. 
Recall (8.1), that is, $\ep_n \leq d_n = 2^{-(N-1)}\,N^{m_0}$. 

If $d_n \geq \frac{1}{2}$, then
\be
H(\ep_n) \leq \log 2 \leq 2 d_n = 2^{-(N-2)}\,N^{m_0}.
\en
Note that
$$
\log\frac{1}{d_n} = m_0 \log\frac{1}{N} + (N-1) \log 2 < N.
$$
Hence, in the other case $d_n \leq \frac{1}{2}$, we have
\be
H(\ep_n) \leq H(d_n) \leq 2d_n \log\frac{1}{d_n} \leq 
2^{-(N - 2)}\, N^{m_0 + 1}.
\en

Comparing (9.10) and (9.11), we see that they can be combined to the
following estimate
$$
H(\ep_n) \leq 2^{-(N - 2)}\, N^{m_0 + 1},
$$
which is valid regardless of whether $d_n$ is greater or smaller than 
$\frac{1}{2}$. 

Using this estimate in (9.9), we finally get
\bee
|\widetilde D_n - D_n| & \leq & 
2^{-(N-2)}\, N^{m_0} \,(4D + 11 + 2N) + 2^{-(N - 2)}\, N^{m_0 + 1} \\
 & = &
2^{-(N-2)}\, N^{m_0} \,(4D + 11 + 3N).
\ene
Since $4D + 11 + 3N < 2^4\, N(D+1)$, we arrive at the desired inequality 
(9.1). 

Thus, Proposition 9.1 is proved.


\vskip5mm
\section{{\bf Integrability of Characteristic Functions $\widetilde f_n$
and their Derivatives}}
\setcounter{equation}{0}

\vskip2mm
Now we turn to the question of quantitative bounds for the modified
characteristic functions $\widetilde f_n$ in terms of the maximal 
entropic distance to normality
$$
D = \max_{k \leq n} D(X_k).
$$

Again, let $X_1,\dots,X_n$ be independent random variables with 
$\E X_k = 0$, $\E X_k^2 = \sigma_k^2$ ($\sigma_k>0$), such that 
$\sigma_1^2 + \dots + \sigma_n^2 = 1$. Moreover, all $X_k$ are assumed
to have absolutely continuous distributions with finite entropies.

We assume that the modified density $\widetilde p_n$ and its 
characteristic function $\widetilde f_n$ have been constructed 
for arbitrary integers $m_0 + 1 \leq N \leq n$.
Put $\sigma = \max_k \sigma_k$.

\vskip5mm
{\bf Proposition 10.1.} {\it If $m_0 \geq 1$ and 
$m_0 + 1 \leq N \leq \frac{1}{2\sigma^2}$, then
\be
\int_{|t| \geq \sqrt{N}}\, |\widetilde f_n(t)|^2\,dt \leq 
C\sqrt{N}\, e^{-cN}
\en
with some positive constants $C$ and $c$, depending on $D$, only. 
}

\vskip5mm
In fact, one can choose the constants to be of the form
$C = e^{2D + 4}$ and $c = c_0 e^{-12\,D}$,
where $c_0$ is a positive absolute factor.

\vskip5mm
{\bf Proof.} Consider any convolution 
$$
\rho = 
(\rho_{10}^{\delta_1} * \rho_{11}^{1-\delta_1}) * \dots * 
(\rho_{N0}^{\delta_N} * \rho_{N1}^{1-\delta_N})
$$ 
participating in the definition of $q_{n0}$, that is, with 
$\delta_1 + \dots + \delta_N > m_0$. It has the Fourier transform
\be
\hat \rho(t)  = \int_{-\infty}^{+\infty} e^{itx}\rho(x)\,dx =
\prod_{k=1}^N \hat \rho_{k0}(t)^{\delta_k}\, \hat \rho_{k1}(t)^{1-\delta_k},
\en
where $\hat \rho_{kj}$ denote the characteristic functions of the
random variables $V_{kj}$ from the median decomposition (4.1) with 
$X = V_k$ ($1 \leq k \leq N$, $j = 0,1$).
In every such convolution there are at least $m_0+1$ terms $\rho_{k0}$ 
for which $\delta_k = 1$. For definiteness, let $k = N$ be one of them, 
so that $\delta_N = 1$. Then, we may write
\be
\hat \rho(t) = 
\hat \rho_{N0}(t) \prod_{k=1}^{N-1} \hat \rho_{k0}(t)^{\delta_k}\, 
\hat \rho_{k1}(t)^{1-\delta_k}.
\en

By Lemma 7.3 $c)$, for all $|t| \geq \sqrt{N}$,
\be
|\hat \rho_{kj}(t)| \leq \exp\big\{-c_0 e^{-12\,D}\big\}
\en
with some absolute constant $c_0>0$. Inserting this in (10.3) and using
$N \geq 2$ leads to
\be
|\hat\rho(t)|^2 \, \leq \, A\,|\hat \rho_{N0}(t)|^2, \qquad
A = \exp\big\{-c_0 e^{-12\,D} N\big\},
\en
where $c_0>0$ is a different absolute constant.

Now, integrate (10.5) over the region $|t| \geq \sqrt{N}$ and 
use Plancherel's formula. Applying the property 
$\rho_{N0}(x) \leq m = m(\rho_N(V_N))$, we get
\be
\int_{|t| \geq \sqrt{N}}\, |\hat \rho(t)|^2\,dt \leq A
\int_{-\infty}^{+\infty} |\hat \rho_{N0}(t)|^2\,dt = 2\pi A\,
\int_{-\infty}^{+\infty} \rho_{N0}(x)^2\,dx \leq 2\pi A\, m.
\en
But, as noted in (7.4), we have $m \leq e^{2D + 2} \sqrt{N}$, so
together with $2\pi < e^2$ (10.6) gives the desired bound
$$
\int_{|t| \geq \sqrt{N}}\, |\hat \rho(t)|^2\,dt \, \leq \, 
e^{2D + 4} \sqrt{N}\, e^{-cN} \qquad 
\big(\,c = c_0 e^{-12\,D}\,\big)
$$ 
for $\hat \rho$. But $\widetilde f_n$ is a finite convex combination of 
such functions, so (10.1) immediately follows.

Thus Proposition 10.1 is proved.

\vskip5mm
Next, we shall extend Propositions 10.1 to the derivatives of 
$\widetilde f_n$, which are needed up to order $s = 3$ in case of 
finite 4-th moments of $X_k$. 
Assume that $s \geq 1$ is an arbitrary integer.

Consider the characteristic functions $\hat \rho$ in (10.2).
Recall that $\widetilde f_n$ represents a convex combination of such
characteristic functions over all sequences 
$\delta = (\delta_1,\dots,\delta_N)$ such that 
$\delta_1 + \dots + \delta_N \geq m_0 + 1$. Hence, it will be sufficient
to derive an estimate, such as (10.1),
for any admissible fixed sequence $\delta$.

Put
$$
u_k = \hat \rho_{k0}^{\delta_k}\, \hat \rho_{k1}^{1-\delta_k} \quad
(1 \leq k \leq N),
$$
which is the characteristic function of the random 
variable $\delta_k V_{k0} + (1-\delta_k)\, V_{k1}$.

Thus, $\hat \rho  = \prod_{k=1}^N u_k$. For the $s$-th derivative 
of the product we write a general polynomial formula
\begin{displaymath}
\hat \rho^{(s)} = 
\sum {s \choose s_1 \dots s_N}\, u_1^{(s_1)} \dots u_N^{(s_N)},
\end{displaymath}
where the summation runs over all integer numbers 
$s_1,\dots, s_N \geq 0$, such that $s_1 + \dots + s_N = s$.

Fix such a sequence $s_1,\dots, s_N$. Note that it contains at most $s$ 
non-zero terms. The sequence $\delta = (\delta_1,\dots,\delta_N)$ 
defining $\rho$ satisfies $\delta_1 + \dots + \delta_N \geq m_0 + 1$. 
Hence, in the row $u_1^{(s_1)}, \dots, u_N^{(s_N)}$ there are at least 
$m_0+1$ terms corresponding to $\delta_k = 1$. Therefore, if $m_0 \geq s$,
there is at least one index, say $k$, for which $\delta_{k} = 1$ and 
in addition $s_k = 0$. For definiteness, let $k = N$, so that
\be
\psi \, \equiv \, 
u_1^{(s_1)} \dots u_N^{(s_N)} \, = \, 
\hat\rho_{N0}\, u_1^{(s_1)} \dots u_{N-1}^{(s_{N-1})}.
\en  

If $s_k>0$, then
$$
|u_k^{(s_k)}(t)| \leq 
\E\, |\delta_k V_{k0} + (1-\delta_k)\, V_{k1}|^{s_k} \leq
\max\{\E\,|V_{k0}|^{s_k}, \E\,|V_{k1}|^{s_k}\}.
$$
But, by the decomposition (7.3) and Jensen's inequality,
$$
\frac{1}{2}\, \E\,|V_{k0}|^{s_k} + \frac{1}{2}\, \E\,|V_{k1}|^{s_k} 
= \E\,|V_k|^{s_k} \leq \E\, |S_n|^{s_k},
$$
so $|u_k^{(s_k)}(t)| \leq 2\,\E\, |S_n|^{s_k}$. Hence, 
\be
\prod_{s_k > 0} |u_k^{(s_k)}(t)| \, \leq \, 
2^s \prod_{s_k > 0} \E\, |S_n|^{s_k} \, \leq \,
2^s \prod_{s_k > 0} \big(\E\, |S_n|^s\big)^{s_k/s} =\, 2^s\,\E\, |S_n|^s.
\en

When $s_k = 0$, we apply the estimate (10.4) on Cramer's constants, 
which may be used in (10.7). Note that (10.4) is fulfilled for at 
least $(N-1) - (s-1) \geq N-m_0$ indices $k \leq N-1$. Hence, 
using also (10.8), we get
$$
|\psi(t)| \leq C\, |\hat\rho_{N0}(t)|\,
\exp\big\{-c_0(N-m_0)\, e^{-12\,D}\big\}, \qquad C = 2^s\,\E\, |S_n|^s.
$$

In case $N \geq 2m_0$, one may simplify this bound by writing
$N-m_0 \geq \frac{N}{2}$. In addition, since the sum of the multinomial
coefficients in the representation of $\hat \rho^{(s)}$ is equal to $N^s$, 
and using Jensen's inequality for the quadratic function, we arrive at
$$
|\hat \rho^{(s)}(t)|^2 \leq \, A\,|\hat\rho_{N0}(t)|^2,\qquad
A = CN^s \exp\big\{-c_0 e^{-12\,D} N\big\},
$$
with some absolute constant $c_0>0$. It remains to integrate
this inequality like in (10.6) over the region $|t| \geq \sqrt{N}$
and apply the estimate (7.4). As a result, we obtain
$$
\int_{|t| \geq \sqrt{N}}\, |\hat \rho^{(s)}(t)|^2\,dt \leq 
Ae^{2D + 4}\,\sqrt{N}.
$$ 

Since $\widetilde f_n$ is a convex combination of the functions
$\hat \rho^{(s)}$, a similar inequality holds for $\widetilde f_n(t)$,
as well. That is,
$$
\int_{|t| \geq \sqrt{N}}\, |\widetilde f_n^{(s)}(t)|^2\,dt
 \ \leq \
2^s\,\E\, |S_n|^s\,e^{2D + 4}\,\exp\big\{-c_0 e^{-12\,D} N\big\}\, N^{s+1/2}.
$$

For $s = 1$ and $s = 2$, we have $\E\, |S_n|^s \leq 1$, while for $s \geq 3$,
one may use Rosenthal's inequality (8.6). In particular, for $s=3$
it gives $\E\, |S_n|^3 \leq 2(1 + L_3)$.

Summarizing the results obtained so far, we have:

\vskip5mm
{\bf Proposition 10.2.} {\it Let 
$m_0 \geq 3$ and $2m_0 \leq N \leq \frac{1}{2\sigma^2}$. Then
\be
\int_{|t| \geq \sqrt{N}}\, |\widetilde f_n^{(s)}(t)|^2\,dt \leq 
C N^{s+1/2}\, e^{-cN} \qquad (s = 1,2)
\en
with positive constants $C$ and $c$, depending on $D$, only.
Moreover, if $L_s$ is finite, $s \geq 3$ integer, and $m_0 \geq s$, then
$$
\int_{|t| \geq \sqrt{N}}\, |\widetilde f_n^{(s)}(t)|^2\,dt \leq 
C \cdot C_s (1+L_s)\, N^{s+1/2}\, e^{-cN}
$$
}

\vskip2mm
Here, the constants $C = e^{2D + 4}$ and $c = c_0 e^{-12\,D}$ 
are of the same form as in Proposition 10.1, and $C_s$ is 
a constant in Rosenthal's inequality (8.6). 
In particular, for $s=3$, we arrive at
\be
\int_{|t| \geq \sqrt{N}}\, |\widetilde f_n'''(t)|^2\,dt \leq 
C (1+L_3)\, N^{7/2}\, e^{-cN}.
\en

Note also that, for $s=0$, (10.9) is true, as well, and returns us
to Proposition 10.1.


\vskip5mm
\section{{\bf Proof of Theorem 1.1 and its Refinement}}
\setcounter{equation}{0}

\vskip2mm
We are now ready to complete the proof of Theorems 1.1-1.2 and
emphasize some of their refinements.
Thus, let $X_1,\dots,X_n$ be independent random variables
with mean zero and finite third absolute moments, having finite
entropies, and such that the sum
$S_n = X_1 + \dots + X_n$ has variance $\Var(S_n) = 1$. 

Our main quantities are the Lyapunov coefficient
$$
L_3 = \sum_{k=1}^n \E\, |X_k|^3
$$
and the maximal entropic distance to normality $D = \max_k D(X_k)$.

To bound the total variation distance $\|F_n - \Phi\|_{{\rm TV}}$
from the distribution $F_n$ of $S_n$ to the standard normal law $\Phi$,
one may apply the general bound (2.1) of Proposition 2.1. 
However, it is only applicable when the characteristic function 
$f_n$ of $S_n$ and its derivative are square integrable. 
But even in the case that, for example, each density $p_n$ of $S_n$ 
is bounded individually, we still could not properly bound
the maximum of the convolutions of these densities explicitly
in terms of $D$ and $L_3$. 
That is why, we are forced to consider modified forms of $p_n$.

Thus, consider these modifications $\widetilde p_n$ together with
their Fourier transforms $\widetilde f_n$ described in Definition 8.1.
By the triangle inequality,
\be
\|F_n - \Phi\|_{{\rm TV}} \leq \|\widetilde F_n - \Phi\|_{{\rm TV}} +
\|\widetilde F_n - F_n\|_{{\rm TV}},
\en
where $\widetilde F_n$ denotes the distribution with density 
$\widetilde p_n$.

In the construction of $\widetilde p_n$ it suffices to take the values
$m_0 = 3$ and $6 \leq N \leq \frac{1}{2\sigma^2}$. Then, by Proposition 8.2, 
\be
\|\widetilde F_n - F_n\|_{{\rm TV}} \, = \, 
\int_{-\infty}^{+\infty} |\widetilde p_n(x) - p_n(x)|\,dx
 \, \leq \,
2^{-(N-2)}\,N^3.
\en
This gives a sufficently good bound on
the last term in (11.1), if $N$ is sufficiently large.

The first term on the right-hand side of (11.1) can be bounded by virtue of 
(2.1), which gives
\be
\|\widetilde F_n - \Phi\|_{{\rm TV}}^2 \, \leq \, \frac{1}{2}\,
\|\widetilde f_n - g\|_2^2 + \frac{1}{2}\,\|(\widetilde f_n)' - g'\|_2^2,
\en
where $g(t) = e^{-t^2/2}$. To estimate the $L^2$-norms, first write
\bee
\frac{1}{2}\,\|\widetilde f_n - g\|_2^2
 & \leq & 
\frac{1}{2} \int_{|t| \leq \sqrt{N}}\,|\widetilde f_n(t) - g(t)|^2\,dt \\
 & & \hskip10mm + \
\int_{|t| > \sqrt{N}}\, |\widetilde f_n(t)|^2\,dt +
\int_{|t| > \sqrt{N}}\, g(t)^2\,dt.
\ene
Since $|\widetilde f_n(t) - f_n(t)| \leq 2^{-(N-2)}\,N^3$, we have
\begin{eqnarray}
\hskip-10mm 
\frac{1}{2} \int_{|t| \leq \sqrt{N}}\,|\widetilde f_n(t) - g(t)|^2\,dt 
 &\leq &
\int_{|t| \leq \sqrt{N}}\,|\widetilde f_n(t) - f_n(t)|^2\,dt +
\int_{|t| \leq \sqrt{N}}\,|f_n(t) - g(t)|^2\,dt \nonumber \\
 &\leq &
 \int_{|t| \leq \sqrt{N}}\,|f_n(t) - g(t)|^2\,dt \, + \, 
2^{-(2N-5)}\,N^{7/2}.
\end{eqnarray}
In addition, by Proposition 10.1,
\be
\int_{|t| \geq \sqrt{N}}\, |\widetilde f_n(t)|^2\,dt \leq 
C\sqrt{N}\, e^{-cN}
\en
with $C = e^{2D + 4}$ and $c = c_0 e^{-12\,D}$,
where $c_0$ is an absolute positive constant.

Using a well-known bound $1 - \Phi(x) \leq \frac{1}{x}\,\varphi(x)$
($x > 0$), we easily get 
$\int_{|t| > \sqrt{N}}\, g(t)^2\,dt < e^{-N}$.
Together with (11.4)-(11.5), and since one may always
assume that $c_0 \leq \frac{1}{2}$, the latter gives
\be
\frac{1}{2}\,\|\widetilde f_n - g\|_2^2 \leq 
\int_{|t| \leq \sqrt{N}}\,|f_n(t) - g(t)|^2\,dt + C \sqrt{N}\, e^{-cN}
\en
with $D$-dependent constants $C = C_0 e^{2D}$ and $c = c_0 e^{-12\,D}$
(where $C_0$ and $c_0$ are numerical).

A similar analysis based on the application of Proposition 8.3 
(cf. Corollary 8.4) and Proposition 10.2 with $s=1$
leads to an analogous estimate 
$$
\frac{1}{2}\,\big\|(\widetilde f_n)' - g'\big\|_2^2 \leq 
\int_{|t| \leq \sqrt{N}}\,|f_n'(t) - g'(t)|^2\,dt + CN^{3/2}\, e^{-cN}.
$$
Together with (11.6) it may be applied in (11.3), and then we get
\bee
\|\widetilde F_n - \Phi\|_{{\rm TV}}^2 
 & \leq & 
\int_{|t| \leq \sqrt{N}}\,|f_n(t) - g(t)|^2\,dt \\
 & & + \ 
\int_{|t| \leq \sqrt{N}}\,|f_n'(t) - g'(t)|^2\,dt + CN^{3/2}\, e^{-cN}.
\ene

It is time to appeal to the classical theorem on the approximation
of $f_n$ by the characteristic function of the standard normal law,
cf. e.g. [R-RR].

\vskip5mm
{\bf Lemma 11.1.} {\it Assume $L_3 \leq 1$. Up to an absolute constant $A$, 
in the interval $|t| \leq L_3^{-1/3}$ we have
$$
|f_n(t) - g(t)| \leq AL_3\,e^{-t^2/4},
$$
and similarly for the first three derivatives of $f_n - g$.
}

\vskip5mm
In fact, the above inequality holds in the larger interval 
$|t| \leq 1/(4L_3)$. But this will not be needed for the present 
formulation of Theorem 1.1.

Thus, if in addition to the original condition 
$6 \leq N \leq \frac{1}{2\sigma^2}$ we require that 
$\sqrt{N} \leq L_3^{-1/3}$, Lemma 11.1 may be applied, and we get
$$
\|\widetilde F_n - \Phi\|_{{\rm TV}} 
 \leq A L_3 + CN^{3/2}\, e^{-cN}.
$$
Using this together with (11.2) in (11.1), we arrive at
\be
\|F_n - \Phi\|_{{\rm TV}} \leq A L_3 + CN^{3/2}\, e^{-cN},
\en
where $A$ is some positive absolute constant, while 
$C = C_0 e^{2D}$ and $c = c_0 e^{-12\,D}$, as before.

\vskip5mm
{\bf Proof of Theorem 1.1.}
To finish the argument, we may take $N = [\frac{1}{2}\,L_3^{-2/3}]$, 
so that $\sqrt{N} \leq L_3^{-1/3}$.
In view of the elementary bound $\sigma \leq L_3^{1/3}$, the condition 
$N \leq \frac{1}{2\sigma^2}$ is fulfilled, as well.
Finally, the condition $N \geq 6$ just restricts us to smaller values 
of $L_3$, and, for example, $L_3 \leq \frac{1}{64}$ would work. Indeed, 
in this case, $\frac{1}{2}\,L_3^{-2/3} \geq 8$, so $N \geq 8$. 

Thus, if $L_3 \leq \frac{1}{64}$, then (11.7) holds true. But since
$N \geq \frac{1}{4}\, L_3^{-2/3}$, the last term in (11.7) is 
dominated by any power of $L_3$ (up to constants). For example,
using $e^x \geq c_1 x^3$ ($x \geq 0$), we get
$$
N^{3/2}\, e^{-cN} \leq \frac{1}{c_1 c^3}\, N^{-3/2}
 \leq \frac{8}{c_1 c^3}\,L_3 = \frac{8}{c_1 c_0^3}\,
e^{36 D}\, L_3.
$$
Hence, (11.7) implies
\be
\|F_n - \Phi\|_{{\rm TV}} \leq C L_3,
\en
with $C = C_0 e^{C_1 D}$, where $C_0, C_1$ are positive numerical
constants.

Finally, if $L_3 > \frac{1}{64}$, (11.8) automatically holds with 
$C = 128$. 

Thus, Theorem 1.1 is proved.

\vskip5mm
Note, however, that the inequality (11.7) contains more information
in comparison with Theorem 1.1. Again assume, as above, that 
$L_3 \leq \frac{1}{64}$ and take $N = [\frac{1}{2}\,L_3^{-2/3}]$. If 
$D \leq \frac{1}{24}\,\log \frac{1}{L_3}$, then
$cN \geq c_0 L_3^{1/2} \cdot \frac{1}{4}\,L_3^{-2/3} = 
c_0' L_3^{-1/6}$ and $C = C_0 e^{2D} \leq C_0 L_3^{-1/12}$.
Hence,
$$
CN^{3/2}\, e^{-cN} \leq C_0 L_3^{-1/12} \cdot L_3^{-1} \cdot
e^{-c_0' L_3^{-1/6}} \leq C_0' L_3
$$
with some absolute constant $C_0'$. As a result, (11.7) yields
$\|F_n - \Phi\|_{{\rm TV}} \leq (A+C_0')\, L_3$. If $L_3 > \frac{1}{64}$, 
(11.8) holds with $C = 128$, and we arrive at:

\vskip5mm
{\bf Theorem 11.2.} {\it Assume that independent random variables $X_k$ 
have mean zero and finite third absolute moments. If they satisfy
$D(X_k) \leq c\,\log \frac{1}{L_3}$ $(1 \leq k \leq n)$,
then
\be
\|F_n - \Phi\|_{{\rm TV}} \leq C L_3,
\en
where $C$ and $c$ are positive absolute constants.
$($One may take $c = \frac{1}{24})$.
}


\vskip5mm
\section{{\bf Proof of Theorem 1.2 and its Refinement}}
\setcounter{equation}{0}

\vskip2mm
In the proof of Theorem 1.2, we apply the general bound (3.1) of 
Proposition 3.1 to the modified densities $\widetilde p_n$
constructed under the same constraints $m_0 = 3$ and 
$6 \leq N \leq \frac{1}{2\sigma^2}$, as in the proof of Theorem 1.1. 
It then gives
$$
\widetilde D_n \, \leq \, \alpha^2 + 
4 \left(\big\|\widetilde f_n - g_\alpha\|_2 + 
\|(\widetilde f_n)''' - g_\alpha'''\big\|_2\right),
$$
where $\widetilde D_n$ is the relative entropy of $\widetilde F_n$
with respect to $\Phi$ and
$$
g_\alpha(t) = g(t)\bigg(1 + \alpha\,\frac{(it)^3}{3!}\bigg), \qquad
\alpha = \sum_{k=1}^n \E X_k^3.
$$

As we know from Proposition 9.1, $\widetilde D_n$ provides a good 
approximation for the entropic distance $D_n = D(S_n)$, namely
$$
|\widetilde D_n - D_n| < 2^{-(N-6)}\, N^4 \,(D + 1).
$$
Hence,
\be
D_n \ \leq \ \alpha^2 + 4 
\left(\big\|\widetilde f_n - g_\alpha\big\|_2 + 
\big\|(\widetilde f_n)''' - g_\alpha'''\big\|_2\right) + 
2^{-(N-6)} N^4 \,(D + 1).
\en

\vskip2mm
On the other hand, the closeness of $f_n$ and $g_\alpha$ on
relatively large intervals is provided by:

\vskip5mm
{\bf Lemma 12.1.} {\it Assume $L_4 \leq 1$. Up to an absolute constant $A$, 
in the interval $|t| \leq L_4^{-1/6}$ we have
\be
|f_n(t) - g_\alpha(t)| \leq AL_4\,e^{-t^2/4},
\en
and similarly for the first four derivatives of $f_n - g_\alpha$.
}

\vskip5mm
Again, we refer to [BR-R], where one can find several variants of such
bounds.

We also use the following elementary relations, cf. e.g. [Pe].

\vskip5mm
{\bf Lemma 12.2.} {\it $\alpha^2 \leq L_3^2 \leq L_4$.
}

\vskip5mm
Now, assume that $L_4 \leq 1$. To estimate the $L^2$-norms in (12.1), 
again write
\begin{eqnarray}
\|\widetilde f_n - g_\alpha\|_2^2
 & \leq & 
\int_{|t| \leq \sqrt{N}}\,|\widetilde f_n(t) - g_\alpha(t)|^2\,dt \nonumber \\
 & & \hskip10mm + \
2 \int_{|t| > \sqrt{N}}\, |\widetilde f_n(t)|^2\,dt +
2 \int_{|t| > \sqrt{N}}\, |g_\alpha(t)|^2\,dt.
\end{eqnarray}

Using $|\widetilde f_n(t) - f_n(t)| \leq 2^{-(N-2)}\,N^3$ and the inequality
(12.2) with $|t| \leq \sqrt{N} \leq L_4^{-1/6}$, we have
\begin{eqnarray}
\hskip-10mm 
\int_{|t| \leq \sqrt{N}}\,|\widetilde f_n(t) - g_\alpha(t)|^2\,dt 
 &\leq &
2 \int_{|t| \leq \sqrt{N}}\,|\widetilde f_n(t) - f_n(t)|^2\,dt +
2 \int_{|t| \leq \sqrt{N}}\,|f_n(t) - g_\alpha(t)|^2\,dt \nonumber \\
 &\leq &
AL_4^2 \, + \, 2^{-(2N-5)}\,N^{7/2}
\end{eqnarray}
with some absolute constant $A$. 

The middle integral on the righ-hand side
of (12.3) has been already estimated in (11.5).

In addition, using $t^6 g(t) \leq 6^3/e^3$, we have
$$
|g_\alpha(t)|^2 = g(t)^2\bigg(1 + \alpha^2\,\frac{t^6}{36}\bigg) <
\big(1 + \alpha^2\big)\, g(t) \leq 2\,g(t),
$$
where we applied Lemma 12.2 together with the assumption $L_4 \leq 1$
(so that $|\alpha| \leq 1$). Hence,
$$
\int_{|t| > \sqrt{N}}\, |g_\alpha(t)|^2\,dt < 
2 \int_{|t| > \sqrt{N}}\, e^{-t^2/2}\,dt <\, 2\,e^{-N/2}.
$$

One may combine this bound with (11.5) and (12.4), and then (12.3) gives
$$
\|\widetilde f_n - g_\alpha\|_2^2 \ \leq \ AL_4^2 + 2^{-(2N-5)}\,N^{7/2}
 + C\sqrt{N}\, e^{-cN} + 4\, e^{-N/2}
$$
with $C = e^{2D + 4}$ and $c = c_0 e^{-12\,D}$ as in (11.5),
where $c_0$ is an absolute positive constant. Since one may always
choose $c_0 \leq \frac{1}{2}$, the above inequality may be simplified as
$$
\|\widetilde f_n - g_\alpha\|_2 \leq AL_4 + CN^{1/4}\, e^{-cN}
$$
with some absolute constant $A$ and $D$-dependent constants 
$C = C_0 e^{2D}$ and $c = c_0 e^{-12\,D}$.

By a similar analysis based on the application of Corollary 8.4 and 
Proposition 10.2 with $s=3$ (cf. (10.10)), we also have an analogous 
estimate 
$$
\|\widetilde f_n''' - g_\alpha'''\|_2 \leq AL_4 + CN^{7/4}\, e^{-cN}.
$$
Hence, (12.1) together with Lemma 12.2 yields
\be
D_n \ \leq \ AL_4 + CN^{7/4}\, e^{-cN},
\en
where $A$ is absolute, and $C = C_0 e^{2D}$ and $c = c_0 e^{-12\,D}$,
as before. The obtained estimate holds true, as long as
$6 \leq N \leq \frac{1}{2\sigma^2}$ and $\sqrt{N} \leq L_4^{-1/6}$
with  $L_4 \leq 1$.

\vskip5mm
{\bf Proof of Theorem 1.2.}
The last condition, $\sqrt{N} \leq L_4^{-1/6}$, is satisfied
for $N = [\frac{1}{2}\,L_4^{-1/3}]$. Then, by the elementary bound 
$\sigma \leq L_4^{1/4}$, we also have $N \leq \frac{1}{2\sigma^2}$. 
The condition $N \geq 6$ restricts us 
to smaller values of $L_4$. If, for example, $L_4 \leq 4^{-6}$,
we have $\frac{1}{2}\,L_4^{-1/3} \geq 8$ and thus $N \geq 8$.

Thus, if $L_4 \leq 4^{-6}$, then (12.5) holds true. But, since
$N \geq \frac{1}{4}\, L_4^{-1/3}$, the last term in (12.5) is 
dominated by any power of $L_4$. In particular, using 
$e^x \geq c_1 x^5$ ($x \geq 0$), we get
$$
N^2\, e^{-cN} \leq \frac{1}{c_1 c^5}\, N^{-3}
 \leq \frac{4^5}{c_1 c^5}\,L_4 = \frac{4^5}{c_1 c_0^5}\,
e^{60\, D}\, L_4.
$$
Hence, (12.5) yields
\be
D_n \leq C L_4
\en
with $C = C_1 e^{2D}\,e^{60\, D} = C_1\,e^{62\, D}$,
where $C_1$ is an absolute constant.

Finally, for $L_4 > 4^{-6}$, one may use the relation $D_n \leq D$ 
(according to the entropy power inequaity), which shows that (12.6) 
holds with $C = 4^6 D$. 

Thus, Theorem 1.2 is proved.

\vskip5mm
Now, again assume, as above, that 
$L_4 \leq 4^{-6}$ and take $N = [\frac{1}{2}\,L_4^{-1/3}]$. If 
$D \leq \frac{1}{48}\,\log \frac{1}{L_4}$, then
$cN \geq c_0 L_4^{1/4} \cdot \frac{1}{4}\,L_4^{-1/3} = 
c_0' L_4^{-1/12}$ and $C = C_0 e^{2D} \leq C_0 L_4^{-1/24}$.
Hence,
$$
CN^{7/4}\, e^{-cN} \leq C_0 L_4^{-1/24} \cdot L_4^{-7/12}
\exp\big\{-c_0' L_4^{-1/12}\big\} \leq C_0' L_4
$$
with some absolute constant $C_0'$. As a result, (12.5) yields
$D_n \leq (A+C_0')\, L_4$. If $L_4 > 4^{-6}$, 
(12.6) holds with $C = 4^6$, and we arrive at another variant
of Theorem 1.2.

\vskip5mm
{\bf Theorem 12.3.} {\it Assume that independent random variables $X_k$ 
have mean zero and finite fourth absolute moments. If they satisfy
$D(X_k) \leq c\,\log \frac{1}{L_4}$ $(1 \leq k \leq n)$,
then
$$
D(S_n) \leq C L_4,
$$
where $C$ and $c$ are certain positive absolute constants.
$($One may take $c = \frac{1}{48})$.
}

\vskip5mm
Let us illustrate this result in the scheme of weighted sums
$$
S_n = a_1 X_1 + \dots + a_n X_n
$$
of independent identically distributed random variables $X_k$, 
such that $\E X_1 = 0$, $\E X_1^2 = 1$, and with coefficients 
such that $a_1^2 + \dots + a_n^2 = 1$. In this case
$L_4 = \E X_1^4 \, \sum_{k=1}^n a_k^4$, so Theorem 12.3 is
applicable, when the last sum is sufficiently small.

\vskip5mm
{\bf Corollary 12.4.} {\it Assume that $X_1$ has density with finite
entropy, and let $\E X_1^4 < +\infty$. If the coefficients satisfy
$$
\sum_{k=1}^n a_k^4 \, \leq \, \frac{1}{\E X_1^4}\,e^{-c D(X_1)},
$$
then
$$
D(S_n) \, \leq \, C\, \E X_1^4 \, \sum_{k=1}^n a_k^4,
$$
where $C$ and $c$ are positive absolute constants.
$($One may take $c = 48)$.
}

\vskip5mm
For example, in case of equal coefficients, so that
$S_n = \frac{X_1 + \dots + X_n}{\sqrt{n}}$, the conclusion becomes
$$
D(S_n) \, \leq \, \frac{C}{n}\, \E X_1^4, \quad {\rm for \ all} \ \ 
n \geq n_1,
$$
which holds true with an absolute constant $C$ and
$n_1 = e^{48 D(X_1)}\,\E X_1^4$.


\vskip5mm
\section{{\bf The Case of Bounded Densities}}
\setcounter{equation}{0}

\vskip2mm
In this Section we give a few remarks about Theorems 1.1-1.2 for the case, 
where the densities of summands $X_k$ are bounded.

First, let us note that, if a random variable $X$ has an absolutely
continuous distribution with a bounded density $p(x) \leq M$,
where $M$ is a constant, and if the variance $\sigma^2 = \Var(X)$
is finite $(\sigma>0)$, then $X$ has finite entropy, and moreover,
\be
D(X) \leq \log\big(M\sigma\sqrt{2\pi e}\big).
\en
Indeed, if $Z$ is a standard normal random variable, and
assuming (without loss of generality) that $\sigma = 1$, we have
$$
D(X) = h(Z)-h(X) = \log\big(\sqrt{2\pi e}\big) +
\int_{-\infty}^{+\infty} p(x)\log p(x)\,dx,
$$
which immediately implies (13.1).

It is wortwile also noticing that, similarly to $D$, the functional
$X \rightarrow M\sigma$ is affine invariant, where 
$M = {\rm ess\,sup}_x\, p(x)$. Therefore, $M\sigma$ does not depend 
neither on the mean or the variance of $X$. In addition, one always has
$M\sigma \geq \frac{1}{\sqrt{12}}$, and the equality is achieved
only for $X$ which is uniformly distributed in a finite interval of 
the real line. (Without proof this lower bound is already mentioned in [St].)

Using (13.1), Theorems 1.1 and 1.2 admit formulations involving
maximum of densities. In the statement below,
let $(X_k)_{1 \leq k \leq n}$ be independent random variables 
with mean zero and variances $\sigma_k^2 = \E X_k^2$ $(\sigma_k > 0$), 
such that $\sum_{k=1}^n \sigma_k^2 = 1$. Let $F_n$
be the distribution function of the sum
$S_n = X_1 + \dots + X_n$.

\vskip5mm
{\bf Corollary 13.1.} {\it Assume that every $X_k$ has density
bounded by $M_k$. If\, $\max_k M_k \sigma_k \leq \widetilde D$, then
\be
\|F_n - \Phi\|_{{\rm TV}} \leq C L_3,
\en
where the constant $C$ depends on $\widetilde D$, only. Moreover,
\be
D(S_n) \leq C L_4.
\en
}

Moreover, one may take $C = C_0 \widetilde D^c$ with
some positive absolute constants $C_0$ and $c$.

In particular, consider the weighted sums
$$
S_n = a_1 X_1 + \dots + a_n X_n
$$
of independent identically distributed random variables $X_k$, 
such that $\E X_1 = 0$, $\E X_1^2 = 1$, and with coefficients 
satisfying $a_1^2 + \dots + a_n^2 = 1$. If $X_1$ has density, 
bounded by $M$, (13.2)-(13.3) yield respectively
$$
\|F_n - \Phi\|_{{\rm TV}} \, \leq \, 
C_M\, \E\, |X_1|^3\, \sum_{k=1}^n |a_k|^3, \qquad
D(S_n) \, \leq \, C_M \, \E X_1^4\, \sum_{k=1}^n a_k^4,
$$
where $C_M$ depends on $M$, only. (One may take $C_M = C_0 M^{c}$).

Moreover, in the i.i.d. case, where
$S_n = \frac{X_1 + \dots + X_n}{\sqrt{n}}$, the last bound may also
be written with an absolute constant $C$, i.e.,
$$
D(S_n) \, \leq \, \frac{C}{n}\, \E X_1^4, \quad {\rm for \ all} \ \ 
n \geq n_1.
$$
Here one may take $n_1 = (M\sqrt{2\pi e})^{48}\, \E X_1^4$.

\vskip5mm
{\bf Acknowledgement.} We would like to thank M. Ledoux for 
pointing us to the relationship between Theorem 1.2 and the 
transportation inequality of E. Rio.

\end{document}